\documentclass[mathpazo]{aamm}
\usepackage{mathtools}
\usepackage{graphicx}
\usepackage{float}
\usepackage{amsmath}
\usepackage{geometry}
\usepackage{mathtools}
\usepackage{bm}
\usepackage{bbm}
\usepackage{caption}
\usepackage{subcaption}
\usepackage{tabularx}
\usepackage{array}
\usepackage{mathrsfs}
\usepackage{bbm}
\usepackage{epstopdf}
\usepackage[colorlinks, linkcolor=red, anchorcolor=green, citecolor=blue]{hyperref}
\usepackage{booktabs}
\usepackage{siunitx}  
\usepackage[compatibility=false]{caption}
\usepackage{algorithm}
\usepackage{algpseudocode}

\algdef{SE}[DOWHILE]{Do}{doWhile}{\algorithmicdo}[1]{\algorithmicwhile\ #1}

\usepackage[format=hang,font=small,textfont=it]{caption}
\numberwithin{equation}{section}

\title{\bfseries The direct-line method  for forward and inverse linear elasticity problems of composite materials in general domains with multiple singularities}

\author[Wei Q H et.~al.]{Qinghua Wei\affil{1},
      Xiaopeng Zhu\affil{1}\comma\corrauth, and Zhongyi Huang\affil{2}}
\address{\affilnum{1}\ College of Science, Liaoning Technical University, Liaoning, China \\
          \affilnum{2}\ Department of Mathematical Sciences, Tsinghua University, Beijing 100084, China}
\emails{{\tt 1303698364@qq.com} (X. Zhu), {\tt 3158563558@qq.com} (Q. Wei),
         {\tt zhongyih@tsinghua.edu.cn} (Z. Huang)}

\date{}
\newtheorem{thm}{Theorem}[section]

\begin{document}

\begin{abstract}
In this work, a combined strategy of domain decomposition and the direct-line method is implemented to solve the forward and inverse linear elasticity problems of composite materials in general domains with multiple singularities.  Domain decomposition technology treats the general domain as the union of some star-shaped subdomains, which can be handled using the direct-line method. The direct-line method demonstrates rapid convergence of the semi-discrete eigenvalues towards the exact eigenvalues of the elliptic operator, thereby naturally capturing the singularities. We also establish optimal error estimates for the proposed method. Especially, our method can handle multiple singular point problems in general regions, which are difficult to deal with by most methods. On the other hand, the inverse elasticity problem is constructed as a energy functional minimization problem with total variational regularization, we use the aforementioned method as a forward solver to reconstruct the lam\'{e} coefficient of multiple singular points in general regions. Our method can simultaneously deduce heterogeneous $\mu$ and $\lambda$ between different materials. Through numerical experiments on three forward and inverse problems, we systematically verified the accuracy and reliability of this method to solve forward and inverse elastic problems in general domains with multiple singularities.


\end{abstract}

\maketitle

\section{Introduction} \label{sec:intro}
With the widespread application of composite materials in various fields, composite elasticity problems can be mathematically formulated as elliptic problems. When there is a sudden change in the geometric shape of interfaces or boundaries, singularities appear in the solutions of these elliptic problems. The stress singularity at these singular points poses a significant challenge to numerical research. The geometric complexity of composite materials, such as the presence of internal cracks and complex interfaces, often leads to elliptic problems with multiple singularities. 

In recent decades, a wide range of methodologies have been proposed for solving problems involving single singular points. Include: Mesh refinement \cite{B-K-P-1979,M-M-C-1989}, Infinite element method \cite{H-1982}, Boundary element method \cite{Li-M-1990,raveendra1991computation}, Incorporation of singular functions in FEM \cite{F-G-W-1973,lin1976finite}, Auxiliary mapping method \cite{oh1995method}, Mixed adaptive FEM \cite{carstensen2000locking}, Galerkin discontinuous adaptive FEM \cite{wihler2006locking}, Modified extended FEM \cite{2016-J}, Finite strip method \cite{C-J-1996}, Stable generalized FEM \cite{gupta2013stable}, Extended isogeometric analysis \cite{ghorashi2012extended}, Extended finite element method(XFEM) \cite{2008-cm}, etc. The methods for solving multiple singular point problems have mainly been developed based on single singular point problems. Up to now, many methods cannot solve the problem of multiple singularities well.

Among all the methods, the direct-line method has received significant attention for solving elliptic problems with single singular point \cite{schiesser2012numerical,xanthis1991method,1997-AMc}. This method does not require prior knowledge of where the singularity is located \cite{2018-JC}, and circumvents the reduced convergence rate in discretization caused by singular points. So far, this method has been used on simple regions with a single singularity \cite{2023-NM,2001-IJOF}. Based on this method, we now develop a framework for handling multiple singularities in more general domain. 

When dealing with elliptic problems that contain multiple singular points, the radial singularity persists in the local polar coordinate system centered at each singular point \cite{2000-JoCP}. To address this challenge, we have proposed a domain decomposition method \cite{2011-CM} suitable for our problem. Specifically, the general domain is decomposed into multiple interconnected subdomains \cite{1990-IM,1987-SJ}, with each subdomain containing exactly one singular point where the direct-line method can be applied individually. By enforcing continuity conditions across the interfaces between adjacent subdomains \cite{2019-IJ,1955-IJFN}, all subdomains are effectively coupled through a global coupling matrix $G$, thereby providing an effective solution strategy for problems with multiple singular points.

The numerical solution of inverse elasticity problems is another cornerstone of our research, with significant applications in non-destructive evaluation for detecting internal flaws such as cracks, voids, and delamination \cite{2005-IP}. Leveraging the fact that full-field displacement measurements $u$ are now accessible via established experimental techniques \cite{grediac2004use, grediac2012full, parker2010imaging}, we develop a numerical framework to simultaneously identify material interfaces and reconstruct the spatially varying Lam\'{e} coefficients. In contrast to existing approaches \cite{2011-IJ,2021-IP}, our method adopts a least-squares formulation \cite{part1} with piecewise constant approximations for the Lam\'{e} parameters. This inverse problem is framed as the minimization of a regularized energy functional that incorporates global variational constraints \cite{vogel2002computational}. To efficiently solve this high-dimensional optimization problem, we employ the gradient-based Adam algorithm \cite{Kingma2015}, which enhances computational efficiency through adaptive moment estimation. It is worth noting that a good forward problem solver is important for solving an inverse problem. Therefore, the combination of direct-line method and domain decomposition technique can be used as a good forward problem solver.

The outline of this paper is as follows: In Section \ref{sec:GDMOL}, we introduce the direct-line coupling with domain decomposition techniques to numerically solve the linear elasticity problems with multiple singular points. Section \ref{sec:inverse} proposes a forward solver based on the method of Section \ref{sec:GDMOL} for solving inverse elasticity problems. In Section \ref{sec:example}, the numerical experiments presented examine three distinct types of regions and validate the efficacy of the proposed method for solving both forward and inverse problems. Finally, We summarize and conclude the feasibility of combining direct-line method with domain decomposition in Section \ref{sec:conclusion}.

\section{The direct-line method for forward problems in general domains with multiple singularities} \label{sec:GDMOL}
\begin{figure}[H]
\centering
\includegraphics[scale=0.65]{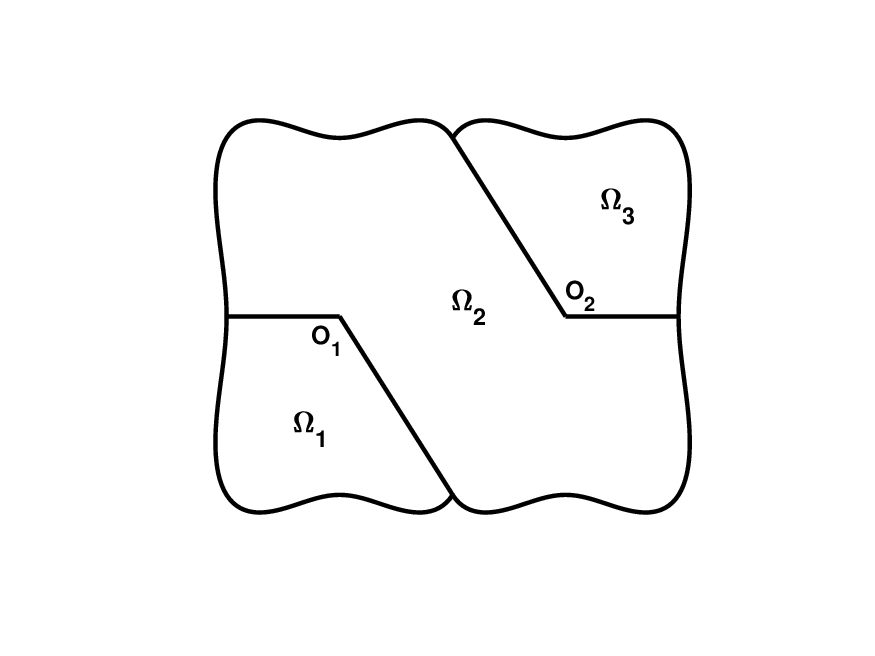}\vspace{-5mm}
\caption{Composite materials in a general domain with multiple singular points and material interfaces}
\label{fig:domain10}
\end{figure}

Consider a general domain $\overline{\Omega} = \bigcup_{k = 1}^K \overline{\Omega}_k \subset R^2$ in a two-dimensional plane, where $\Omega$ represents a domain composed of $K$ different materials and $\Omega_k$ represents the $k$-th material. Let $O_{\tilde{k}}$, $ \tilde{k} =1,...,{\tilde{K}}$ denote several singular points. The boundary $\Gamma = \partial \Omega$ can be approximated as star-shaped relative to some singular points $O_{\tilde{k}}$. Considering the inherent complexity of arbitrary domains, in this work we hypothesize that all material interfaces coincide at singular points $O_{\tilde{k}}$, then each interface $\partial\overline{\Omega}_{k} \bigcap \partial\overline{\Omega}_{k+1}$ is a simple polyline $L_{k}$. see Figure \ref{fig:domain10}. 
The equilibrium state of composite materials under linear elasticity assumptions is governed by Navier's equations defined on a geometrically general domain $\Omega$, subject to Dirichlet boundary constraints:
\begin{equation} \label{composite-eq1-domain}
\begin{aligned}
-\nabla\cdot \sigma_k & = 0, \quad \text{in } \Omega_k,\\
u & = f, \quad \text{on } \Gamma,\\
u_{k} & = u_{k+1}, \quad \text{on } L_{k},\\
\sigma_{k} \cdot {n}_k & = \sigma_{k+1} \cdot {n}_k \quad \text{on } L_{k},
\end{aligned}
\end{equation}
with
$$
\sigma_k = 2 \mu_k \varepsilon(u_k) + \lambda_k \nabla \cdot u_k I, \quad
\varepsilon(u_k) = \frac{1}{2} \left( \nabla u_k + (\nabla u_k)^T \right),k = 1, \ldots, K
$$
where given a vector-valued function $f$ on $\Gamma$, $\sigma|_{\Omega_k} = \sigma_k$ and $\sigma$ is the stress tensor, $u|_{\Omega_k} = u_k$ and $u$ is the displacement, ${n}_k$ is the unit normal vector of $L_k$. The Lam\'e coefficients $\mu$ and $\lambda$ are defined as piecewise-constant functions and taking distinct positive values $(\mu_k, \lambda_k)$ in each subdomain $\Omega_k$.

The forward problem associated with \eqref{composite-eq1-domain} involves computing the displacement field $u$ for known the Lam{\'e coefficients $(\mu,\lambda)$. On the other hand, the inverse linear elasticity problem consists in reconstructing the Lam{\'e coefficients $(\mu,\lambda)$ from the observed displacement data $u$.
We assume that the entire domain $\Omega$ can be separated into $\tilde{K}$ subdomains $\Omega^{(\tilde{k})}$ such that each subdomain contains only one singular point. Let the boundary $\Gamma^{(\tilde{k})} = \partial \Omega^{(\tilde{k})}$ be star-shaped w.r.t. the singular point $O_{\tilde{k}}$ and can be parameterized as a (piecewise) $C^1$ function of the angular variable $\phi^{(\tilde{k})}$, this can be represented as $\tilde{r}^{(\tilde{k})}(\phi^{(\tilde{k})})$, obviously $\tilde{r}^{(\tilde{k})}(0) = \tilde{r}^{(\tilde{k})}(2\pi)$ and $\tilde{r}^{(\tilde{k})}(\phi^{(\tilde{k})}) \geq r_0>0$ for any discretization of angular variables $0 \leq \phi^{(\tilde{k})} \leq 2\pi$. Then we could introduce different curvilinear coordinates in each $\Omega^{(\tilde{k})}$, which constitutes the first step of the direct-line method.
\begin{equation}
(x^{(\tilde{k})}, y^{(\tilde{k})}) = e^{\rho^{(\tilde{k})}}\tilde{r}^{(\tilde{k})}(\phi^{(\tilde{k})}) ( \cos(\phi^{(\tilde{k})}), \sin(\phi^{(\tilde{k})})), 0 \leq \phi^{(\tilde{k})} < 2\pi, -\infty < \rho^{(\tilde{k})} \leq 0.
\label{coordinate}
\end{equation}
\subsection{Domain decomposition and coupling of two singular points}
\label{2.1}
In this subsection, we present a coupled direct-line method with domain decomposition to numerically solve linear elastic problems of composite materials possessing two singular points.  
\begin{figure}[H]
\centering
\includegraphics[scale=0.65]{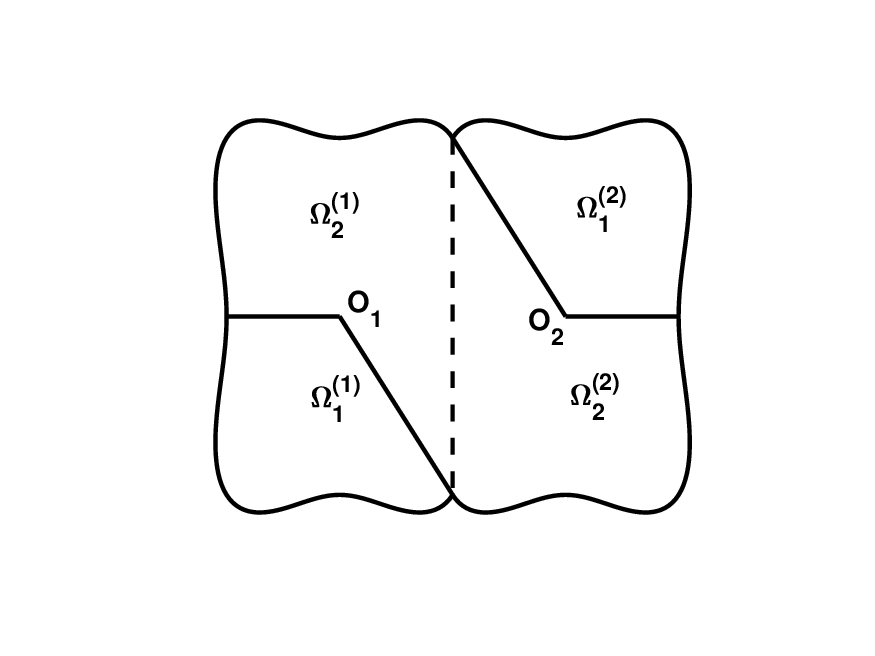}\vspace{-5mm}
\caption{Decomposition of domain $\Omega$ into two subdomains, each containing one singular point}
\label{fig:domain1}
\end{figure}
Assume that the entire domain $\Omega$ can be separated into two subdomains $\Omega^{(1)}$ and $\Omega^{(2)}$ such that each subdomain contains only one singular point. The dashed line $L_c = \partial \Omega^{(1)} \bigcap \partial \Omega^{(2)}= \partial \Omega_2^{(1)} \bigcap \partial\Omega_2^{(2)}$, $\Omega^{(\tilde{k})}=\bigcup \Omega^{(\tilde{k})}_i,i=1,2$. And the polyline $L_{\tilde{k}}=\partial\Omega_1^{\tilde{k}}\bigcap\partial\Omega_2^{\tilde{k}}$. This completes our domain decomposition (see Figure \ref{fig:domain1}).
Stress singularities will still occur at each singular point and the singularities exhibit radial distribution characteristics. Following the domain decomposition principle, problem \eqref{composite-eq1-domain} is equivalent to the following coupled PDE systems:
\begin{equation}
\begin{aligned}
-\nabla\cdot \sigma^{(1)} & = 0, \quad \text{in } \Omega^{(1)}, \\
-\nabla\cdot \sigma^{(2)} & = 0, \quad \text{in } \Omega^{(2)}, \\
u^{(1)} & = f, \quad \text{on } \partial \Omega^{(1)}\backslash\partial\Omega^{(2)},\\
u^{(2)} & = f, \quad \text{on } \partial \Omega^{(2)}\backslash\partial\Omega^{(1)},\\
u^{(1)} & = u^{(2)},\quad \text{on } L_c,\\
\sigma^{(1)}\cdot {n}_c & = \sigma^{(2)}\cdot {n}_c,\quad \text{on } L_c,\\
u^{(\tilde{k})}_1 & = u^{( \tilde{k})}_2, \quad \text{on } L_{\tilde{k}}& \\
\sigma^{(\tilde{k})}_1\cdot {n}_{\tilde{k}} & = \sigma^{(\tilde{k})}_2\cdot {n}_{\tilde{k}} \quad \text{on } L_{\tilde{k}}&
\end{aligned}
\label{composite-eqpde}
\end{equation}

We now proceed with the numerical solution of problem \eqref{composite-eqpde}. It is crucial to note that the overlapping boundary portions between domains $\Omega^{(1)}$ and $\Omega^{(2)}$ serve as the key coupling mechanism to connect the subdomains. In subdomain $\Omega^{(1)}$, we discretize the interval $[0, 2\pi]$ as $0 = \phi_1^{(1)} < \phi_2^{(1)} < \cdots < \phi_{M+1}^{(1)} = 2\pi$, denote $\Phi=\left\{1,...,M\right\}$ and let $I^{(1)}=\left\{j_1,...,j_{m}\right\}\subset\Phi^{(1)}$ be those such that $\forall j_i \in I^{(1)},(0,\phi_{j_i}^{(1)}) \in \partial \Omega^{(1)}\cap\partial\Omega^{(2)}$ in the curvilinear coordinate system \eqref{coordinate} of $\Omega^{(1)}$. Similarly, let $I^{(2)}=\left\{k_1,...,k_{m}\right\}\subset\Phi^{(2)}$ be those such that $\forall k_i \in I^{(2)},(0,\phi_{k_i}^{(2)}) \in  \partial\Omega^{(1)} \cap \partial \Omega^{(2)}$ in the curvilinear coordinate system \eqref{coordinate} of $\Omega^{(2)}$. 
The numerical solution in the subdomains $\Omega^{(\tilde{k})}$ can be expressed as
\begin{equation} \label{composite-eq2222}
\begin{aligned}
u_h^{(\tilde{k})}{(\rho,\phi) = N^{(\tilde{k})}{(\phi)^{T}}S^{(\tilde{k})}{(\rho)}\alpha^{(\tilde{k})}, \tilde{k}=1,2}.
\end{aligned}
\end{equation}

The construction of $N^{(\tilde{k})},S^{(\tilde{k})},\alpha^{(\tilde{k})},\tilde{k}=1,2$ follows a consistent formulation pattern in reference \cite{2023-NM}. We similarly obtain $N^{(\tilde{k})}$ and $S^{(\tilde{k})}$. The coefficient vector $\alpha^{(\tilde{k})}=[\alpha_1^{(\tilde{k})}, \alpha_2^{(\tilde{k})},...,\alpha_{2M}^{(\tilde{k})}]^T$ is determined by enforcing the boundary conditions and the interface conditions that couple the subdomains. A coupling matrix of dimension $4M\times 4M$ is constructed to enforce these constraints.
\begin{equation} \label{composite-eq11}
\begin{aligned}
G = 
\begin{bmatrix}
V_{\Phi\backslash{I^{(1)}}}^{(1)} && 0
\\
0 && V_{\Phi\backslash{I^{(2)}}}^{(2)}
\\
V_{{I^{(1)}}}^{(1)} && -V_{{I^{(2)}}}^{(2)}
\\
H_{{I^{(1)}}}^{(1)} && H_{{I^{(2)}}}^{(2)}
\end{bmatrix}
\end{aligned}
\end{equation}
Where ${\Phi\backslash{I}^{(1)}}=\left\{n_1,...,n_{M-m}\right\}$ represents the remaining portion of the boundary of $\Omega^{(1)}$ excluding $I^{(1)}$. Similarly, let ${\Phi\backslash{I}^{(2)}}=\left\{l_1,...,l_{M-m}\right\}$. Let $V_{{I^{(1)}}}^{(1)}$ be rows of $S^{(1)}$ whose indices are $\left\{j_1,...,j_{m}\right\}$, $V_{{I^{(2)}}}^{(2)}$ be rows of $S^{(2)}$ whose indices are $\left\{k_1,...,k_{m}\right\}$, $V_{\Phi\backslash{I^{(1)}}}^{(1)}$ be rows of $S^{(1)}$ whose indices are $\left\{n_1,...,n_{M-m}\right\}$, $V_{\Phi\backslash{I^{(2)}}}^{(2)}$ be rows of $S^{(2)}$ whose indices are $\left\{l_1,...,l_{M-m}\right\}$. Let $H_{I^{(1)}}^{(1)}$ be rows of $H^{(1)}$ whose indices are $\left\{j_1,...,j_{m}\right\}$, $H_{{I^{(2)}}}^{(2)}$ be rows of $H^{(2)}$ whose indices are $\left \{k_1,...,k_{m}\right\}$, where
$$
H^{\tilde{(k)}} = \sum_{\tilde{k}=1}^{{\tilde{K}}} \int_{\theta_{\tilde{k}}}^{\theta_{\tilde{k}+1}} \left(   N(\phi) \Psi_1 N'(\phi)^T  + N(\phi) \Psi_2 N(\phi)^T S'(0)S(0)^{-1} \right)^{\tilde{(k)}} d\phi,
$$
from reference \cite{2023-NM}. Set the $4M\times 1$ vector
\begin{equation} 
\label{composite-eq00}
F=[F^{(1)T}_{\Phi\backslash{I^{(1)}}},F^{(2)T}_{\Phi\backslash{I^{(2)}}},0...,0]^{T},
\end{equation}
where $F_{\Phi\backslash{I^{(1)}}}^{(1)}$ be rows of $F^{(1)}$ whose indices are $\left\{n_1,...,n_{M-m}\right\}$, $F_{\Phi\backslash{I^{(2)}}}^{(2)}$ be rows of $F^{(2)}$ whose indices are $\left\{l_1,...,l_{M-m}\right\}$, $F^{(\tilde{k})}$ is the Dirichlet boundary condition in $\Omega^{(\tilde{k})}$. Then by the boundary conditions and interface conditions, $\alpha^{(1)}$ and $\alpha^{(2)}$ satisfy
\begin{equation} 
\label{composite-eq33}
G[\alpha^{(1)T},\alpha^{(2)T}]^{T}=F
\end{equation}
After solving \eqref{composite-eq33}, we would obtain the numerical solution $u^{(1)}_h(\rho,\phi)$ and $u^{(2)}_h(\rho,\phi)$ in subdomain $\Omega^{(1)}$ and $\Omega^{(2)}$.
\subsection{Domain decomposition and coupling of four singular points}
\label{2.2}
In this subsection, we extend the methodology for elliptic problems with two singular points to the more challenging case of four singular points. We consider the following problem \eqref{composite-eq1-domain}. We assume that the entire domain $\Omega=\bigcup \Omega_k,k=1,...,5$ (see Figure \ref{fig:domain2a}) can be partitioned into four subdomains $\Omega=\bigcup\Omega^{(\tilde{k})},\tilde{k}=1,..,4$, and $\Omega^{( \tilde{k})}=\bigcup \Omega^{(\tilde{k})}_i,i=1,2$ (see Figure \ref{fig:domain2b}).

\begin{figure}[H]
\centering
\begin{subfigure}[b]{0.46\textwidth}  
\includegraphics[width=\textwidth]{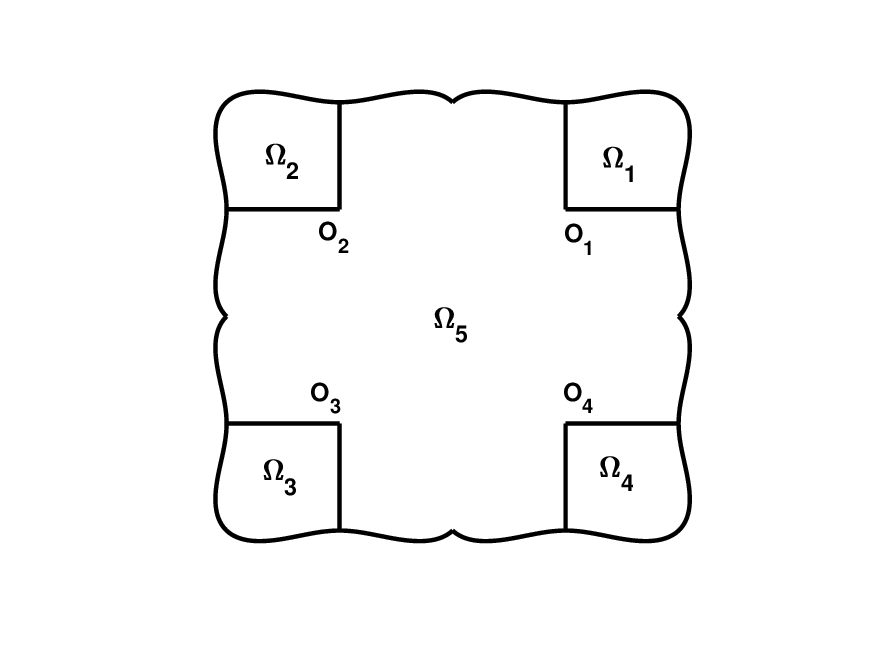}  
\caption{Domain $\Omega$}  
\label{fig:domain2a}
\end{subfigure}
\begin{subfigure}[b]{0.46\textwidth}  
\includegraphics[width=\textwidth]{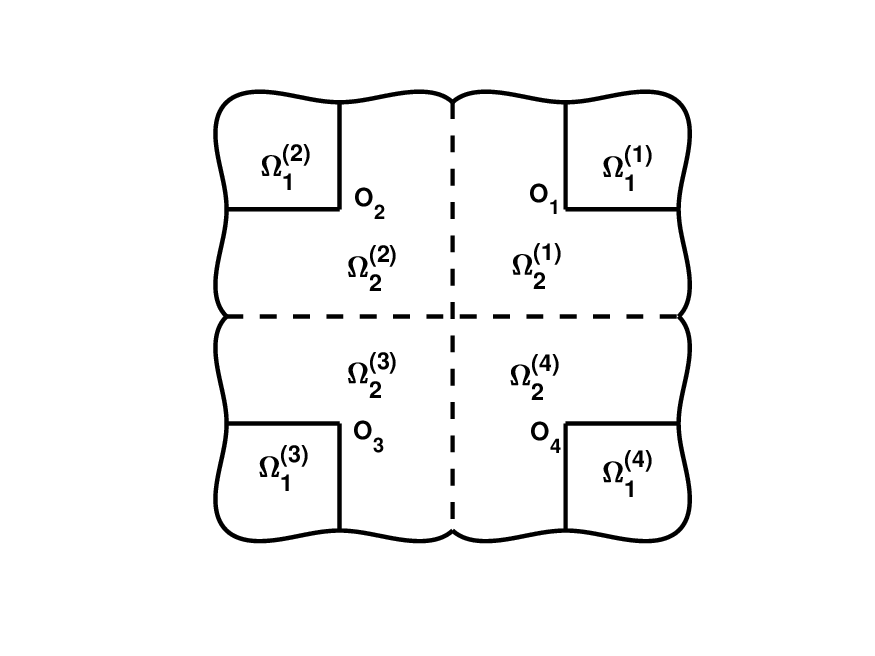}  
\caption{Subdomain $\Omega^{(\tilde{k})}$}  
\label{fig:domain2b}
\end{subfigure}
\caption{(a) Original domain with four singular points; (b) Partition into four subdomains}
\end{figure}
Following the domain decomposition principle, problem \eqref{composite-eq1-domain} is equivalent to four coupled sub-problems. Under these conditions, \eqref{composite-eq1-domain} becomes equivalent to the following coupled PDE system
\begin{equation}
\begin{aligned}
-\nabla\cdot \sigma^{(\tilde{k})} & = 0, \quad \text{in } \Omega^{(\tilde{k})}, \\
u^{(\tilde{k})} & = f, \quad \text{on } \partial \Omega^ {(\tilde{k})} \cap \partial \Omega,\\
u^{(1)} & = u^{(2)},\quad \text{on } L_{c1},\\
u^{(2)} & = u^{(3)},\quad \text{on } L_{c2},\\
u^{(3)} & = u^{(4)},\quad \text{on } L_{c3},\\
u^{(4)} & = u^{(1)},\quad \text{on } L_{c4},\\
\sigma^{(1)}\cdot {n}_{c1} & = \sigma^{(2)}\cdot {n}_{c1},\quad \text{on } L_{c1},\\
\sigma^{(2)}\cdot {n}_{c2} & = \sigma^{(3)}\cdot {n}_{c2},\quad \text{on } L_{c2},\\
\sigma^{(3)}\cdot {n}_{c3} & = \sigma^{(4)}\cdot {n}_{c3},\quad \text{on } L_{c3},\\
\sigma^{(4)}\cdot {n}_{c4} & = \sigma^{(1)}\cdot {n}_{c4},\quad \text{on } L_{c4},\\
u^{(\tilde{k})}_1 & = u^{( \tilde{k})}_2, \quad \text{on } L_{\tilde{k}} \\
\sigma^{(\tilde{k})}_1\cdot {n}_{\tilde{k}} & = \sigma^{(\tilde{k})}_2\cdot {n}_{\tilde{k}} \quad \text{on } L_{\tilde{k}}
\label{eq-2.pde}
\end{aligned}
\end{equation}
where the dashed line $L_c=\bigcup L_{ci},i=1,...,4$, represent the interfaces between subdomains, $L_{c1}=\partial \Omega^{(1)}\bigcap\partial\Omega^{(2)}$, $L_{c2}=\partial \Omega^{(2)}\bigcap\partial\Omega^{(3)}$, $L_{c3}=\partial \Omega^{(3)}\bigcap\partial\Omega^{(4)}
$, $L_{c4}=\partial \Omega^{(4)}\bigcap\partial\Omega^{(1)}$, and the polyline $L_{\tilde{k}}=\partial\Omega^{(\tilde{k})}_1\bigcap\partial\Omega^{(\tilde{k})}_2$.

We now proceed with the numerical solution of problem \eqref{eq-2.pde}. In subdomain $\Omega^{(1)}$, we discretize the interval $[0, 2\pi]$ as $0 = \phi_1^{(1)} < \phi_2^{(1)} < \cdots < \phi_{M+1}^{(1)} = 2\pi$, denote $\Phi=\left\{1,...,M\right\}$ and let $I^{(1,2)}=\left\{j_1,...,j_{m}\right\}\subset\Phi$ be those such that $\forall j_i \in I^{(1,2)},(0,\phi_{j_i}^{(1)}) \in \partial \Omega^{(1)}\cap\partial\Omega^{(2)}$ in the curvilinear coordinate system \eqref{coordinate} of $\Omega^{(1)}$, and let $I^{(1,4)}= \left \{j'_1,...,j'_{m}\right\} \subset \Phi$ be those such that $\forall j'_i \in I^{(1,4)},(0,\phi_{j'_i}^{(1)}) \in \partial \Omega^{(1)}\cap\partial\Omega^{(4)}$ in the curvilinear coordinate system \eqref{coordinate} of $\Omega^{(1)}$. We will not go into details about the other subdomains.
We now proceed with the numerical solution of equation \eqref{eq-2.pde}. The solutions can also be expressed as
\begin{equation} 
\label{composite-eq1111}
u_h^{(\tilde{k})}{(\rho,\phi) = N^{(\tilde{k})}{(\phi)^{T}}S^{(\tilde{k})}{(\rho)}\alpha^{(\tilde{k})}, \tilde{k}=1,..,4}.
\end{equation}

The construction of $N^{(\tilde{k})},S^{(\tilde{k})},\alpha^{(\tilde{k})},\tilde{k}=1,...,4$ follow a consistent formulation pattern. We similarly obtain $N^{(\tilde{k})}$ and $S^{(\tilde{k})}$. To obtain $\alpha^{(\tilde{k})}=[\alpha_1^{(\tilde{k})},\alpha_2^{(\tilde{k})},...,\alpha_{2M}^{(\tilde{k})}]$, we incorporate both boundary conditions and interface conditions that govern the subdomain interactions. A coupling matrix of dimension $8M\times 8M$ is constructed to enforce these constraints.

\begin{equation} \label{composite-eq1}
\begin{aligned}
G = 
\begin{bmatrix}
V_{\Phi\backslash{I_1}}^{(1)} && 0 && 0 && 0\\
0 && V_{\Phi\backslash{I_2}}^{(2)} && 0 && 0\\
0 && 0 && V_{\Phi\backslash{I_3}}^{(3)} && 0\\
0 && 0 && 0 && V_{\Phi\backslash{I_4}}^{(4)}\\
V_{I^{(1,2)}}^{(1)} && -V_{I^{(2,1)}}^{(2)} && 0 && 0\\
0 && V_{I^{(2,3)}}^{(2)} && -V_{I^{(3,2)}}^{(3)} && 0\\
0 && 0 && V_{I^{(3,4)}}^{(3)} && -V_{I^{(4,3)}}^{(4)}\\
-V_{I^{(1,4)}}^{(1)} && 0 && 0 && V_{I^{(4,1)}}^{(4)}\\
H_{I^{(1,2)}}^{(1)} && H_{I^{(2,1)}}^{(2)} && 0 && 0\\
0 && H_{I^{(2,3)}}^{(2)} && H_{I^{(3,2)}}^{(3)} && 0\\
0 && 0 && H_{I^{(3,4)}}^{(3)} && H_{I^{(4,3)}}^{(4)}\\
H_{I^{(1,4)}}^{(1)} && 0 && 0 && H_{I^{(4,1)}}^{(4)}
\end{bmatrix}
\end{aligned}
\end{equation}
where $I_1=I^{(1,2)}\bigcup I^{(1,4)}$, $I_2=I^{(2,3)}\bigcup I^{(2,1)}$, $I_3=I^{(3,4)}\bigcup I^{(3,2)}$, $I_4=I^{(4,1)}\bigcup I^{(4,3)}$. ${\Phi \backslash{I}_1} = \left \{n_1,...,n_{M-2m}\right\}$ represents the remaining portion of the boundary of $\Omega^{(1)}$ excluding $I_1$. Here, matrix $G$ is not elaborated upon further, similar to Subsection \ref{2.1}.\\
and the $8M$-by-1 vector
\begin{equation} \label{composite-eq11111111111}
\begin{aligned}
F=[F^{(1)T}_{\Phi\backslash{I_1}},F^{(2)T}_{\Phi\backslash{I_2}},F^{(3)T}_{\Phi\backslash{I_3}},F^{(4)T}_{\Phi\backslash{I_4}},0,...,0]^{T}
\end{aligned}
\end{equation}
Then, by the boundary conditions and interface conditions, $\alpha^{(\tilde{k})}$ satisfies
\begin{equation} \label{composite-eq100}
\begin{aligned}
G[\alpha^{(1)T},\alpha^{(2)T},\alpha^{(3)T},\alpha^{(4)T}]^{T}=F
\end{aligned}
\end{equation}
After solving \eqref{composite-eq100}, we would obtain the numerical solution $u^{(\tilde{k})}_h(\rho,\phi)$ and in subdomain $\Omega^{(\tilde{k})}$.
In subsections \ref{2.1} and \ref{2.2}, the specific form of $G$ varies with the geometry of the domain.
\subsection{Error estimation}
\begin{thm}
\label{thm-error}
Assume the use of linear elements, there exists a constant $C > 0$ independent of the mesh size $h$ such that the following error estimate holds:
\end{thm}
\begin{equation} \label{eq: error estimate}
\|u-u_h\|^2_*\leq Ch^2\sum_{j=1}^2\iint_{\Omega} \left( |\nabla u_j|^2+\sum_{|\beta|=2}(x^2+y^2)|D^{\beta}u_j|^2 \right) dxdy,
\end{equation}
where
$|\beta| = \beta_1 + \beta_2$, $\beta = (\beta_1, \beta_2)^T \in {N}^2$, $D^{\beta} u_j = \frac{\partial^{|\beta|} u_j}{\partial x^{\beta_1} \partial y^{\beta_2}}$ and for any $w, v \in H^1(\Omega) \times H^1(\Omega)$,
\begin{gather}
\|w\|_* = E(w,w)^{1/2} + |\psi_1(w)| + |\psi_2(w)| + |\psi_3(w)|, \notag\\
\psi_1(w) = \int_{\Gamma}w_1ds,\quad
\psi_2(w) = \int_{\Gamma}w_2ds,\quad
\psi_3(w) = \int_{\Gamma}(w_1y-w_2x)ds,\notag
\end{gather}
$$
\begin{aligned}
E(w,v) = \sum_{k=1}^K \iint_{\Omega_k} & \Big{(} \lambda_k (\nabla \cdot w^k) \cdot (\nabla \cdot v^k) + 2\mu_k \left( \partial_x w_1^k \partial_x v_1^k + \partial_y w_2^k \partial_y v_2^k \right) \\
& + \mu_k \left( \partial_y w_1^k + \partial_x w_2^k \right) \left( \partial_y v_1^k + \partial_x v_2^k \right) \Big{)}dxdy.
\end{aligned}
$$

Recall that $u$ belongs to $H^1(\Omega)$, but the solution of $u$ is not in the $H^2(\Omega)$ space, so $|\beta| = 2$ is integrable, and when $|\beta| = 2$, $|D^{\beta} u|^2$ is not integrable. However, considering that $u$ has only singularities of the form $\mathcal{O}(r^{b_j})$ in the $r$ direction near singular points, where $0 < \textbf{Re}(b_j) < 1$ and $\phi$ is smooth, it follows that $r^2 |D^{\beta} u|^2$ is integrable for $|\beta| = 2$. In Theorem \ref{thm-error}, the integral \eqref{eq: error estimate} is finite. Therefore, from Theorem \ref{thm-error}, it can be concluded that even if $u \notin H^2(\Omega)$, we can still use linear elements for semi-discretization to obtain optimal first-order convergence in the energy norm $|| \cdot ||_{*}$. Furthermore, we assume that for each region there are normal numbers $(\lambda,\mu)$, we can conclude $|| \sigma - \sigma^h ||_{L^2} \le C || u - u^h ||_{*}$, where $h$ is independent of the constant $C$, and $\sigma^h$ denotes the computed stress tensor associated with the discrete displacement field $u^h$. At the same time, we can use linear elements to achieve first-order convergence in $L^2$ norm for singular stress $\sigma$. 
\section{Inverse Elasticity Problems} \label{sec:inverse}
The inverse problem aims to reconstruct the Lam\'{e} coefficients from displacement field measurement data. The study considers composite elastic materials within a general domain, where material interfaces are represented as polygonal lines connecting singular points to boundary points. The main goal is to restore the values of the Lam\'{e} parameters and the geometric shapes of these interfaces simultaneously on all subdomains.
\subsection{A regularized minimization problem}
We let $u[\ell]$ be a solution to \eqref{composite-eq1-domain}, here define the Lam\'{e} coefficients $\ell = (\mu, \lambda)$. For problems involving abrupt material property transitions, as manifested through discontinuous Lam\'{e} coefficients, We are given certain constants $C_1, C_2 > 0$, and for any $v \in L^1(\Omega)$, which is defined as
\begin{equation}
\ell \in \Lambda := \left\{ \ell = (\mu, \lambda) \in L^{\infty}(\Omega) \times L^{\infty}(\Omega)| C_1 \le \mu, \lambda \le C_2, \text{TV}(\mu) < \infty, \text{TV}(\lambda) < \infty \right\},
\end{equation}

\begin{equation}
\text{TV}(v) = \sup \left\{ \iint_{\Omega} \left( v(\nabla \cdot w) \right) dxdy \Big{|} w \in (C_0^1(\Omega))^2, |w(x)| \leq 1, \forall x \in \Omega \right\},
\end{equation}
where $\text{TV}(v)$ is the total variation of $v$ \cite{vogel2002computational} and $|\cdot|$ denotes the Euclidean norm of a vector. In the inverse elasticity problem with full-field data, we consider a noisy measurement $z(x, y)$ of the displacement field $u[\ell^{\star}](x, y)$, where $(x, y) \in \Omega$ and $\ell^{\star} = (\mu^{\star}, \lambda^{\star}) \in \Lambda$. The objective is to reconstruct the unknown parameters $\ell^{\star}$ from the measured data $z$.
Construct an optimization problem to solve the reconstruction of $\ell^{\star}$, where the goal is to minimize a regularized energy functional \cite{ident-lame}. Since the coefficients of each subdomain are not identical, that is, they have jump discontinuities in the Lam\'{e} coefficients. Thus, by using the total variation method, we can effectively recover the Lam\'{e} coefficients \cite{vogel2002computational}. We get an optimization problem
\begin{equation} \label{minpro}
\min_{\ell \in \Lambda} J(\ell) := \frac{1}{2} \iint_{\Omega} \left( \lambda|\nabla \cdot (u[\ell]-z)|^2 + 2 \mu |\varepsilon(u[\ell]-z)|^2 \right) dxdy + \eta \left( \text{TV}(\mu) + \text{TV}(\lambda) \right),
\end{equation}
where $\eta$ is the regularization parameter and greater than 0, and let
$$
J_0(\ell) = \frac{1}{2} \iint_{\Omega} \left( \lambda|\nabla \cdot (u[\ell]-z)|^2 + 2\mu|\varepsilon(u[\ell]-z)|^2  \right) dxdy.
$$
A theorem in the literature \cite{part1} gives the existence of solutions for \eqref{minpro}. Here is a brief explanation.
\begin{thm}\label{theory 1}
The minimization problem \eqref{minpro} exists at least one solution $\ell^0\in \Lambda$.
\end{thm}

\subsection{Numerical program}
Continuing our analysis within the curvilinear coordinate established in \eqref{coordinate}. In the subdomain $\Omega^{\tilde{(k)}}$, for a fixed $\rho^{\tilde{(k)}}$, they are piecewise constant functions along the $\phi^{\tilde{(k)}}$-direction, and for a fixed $\phi^{\tilde{(k)}}$, both $\mu^{\star}$ and $\lambda^{\star}$ remain constant along the $\rho^{\tilde{(k)}}$-direction. Hence we let $0 = \varphi_0^{\tilde{(k)}} < \varphi_1^{\tilde{(k)}} < \ldots < \varphi_{m'}^{\tilde{(k)}} = 2\pi$ and define the discrete parameter space
\begin{equation} \label{discrete admissible set}
\Lambda_{h'} = \left\{ \ell_{h'}(\rho, \phi) = (\mu_{h'}, \lambda_{h'}) \in \Lambda \Bigg{|}
\begin{aligned}
& \mu_{h'}, \lambda_{h'} \text{ are constants almost everywhere} \\
& \text{on } (-\infty,0] \times [\varphi_{j-1}^{\tilde{(k)}}, \varphi_{j}^{\tilde{(k)}}) \text{ for } j = 1,\ldots,m'.
\end{aligned} 
\right\},
\end{equation}
where $h' = \max\limits_{1 \leq \tilde{k}< \tilde{K}}\max \limits_{1 \leq j< m} |\varphi_{j}^{\tilde{(k)}} - \varphi_{j-1}^{\tilde{(k)}}|$. We discretize \eqref{minpro} by the following discrete minimization problem:
\begin{equation} \label{num-minpro}
\begin{aligned}
\min_{\ell_{h'} \in \Lambda_{h'}} J(\ell_{h'}) & = \frac{1}{2} \sum_{\tilde{k}=1}^{\tilde{K}}\sum_{j = 1}^{m'} \int_{-\infty}^0 \int_{\varphi_{j-1 }^{\tilde{(k)}}}^{\varphi_{j}^{\tilde{(k)}}} \left( 2 \mu_{h'} |\varepsilon(u[\ell_{h'}]-z)|^2 + \lambda_{h'} |\nabla \cdot (u[\ell_{h'}]-z)|^2 \right)\\
&e^{2\rho^{\tilde{(k)}}} \tilde{r}^{\tilde{(k)}}(\phi)^2 d\phi^{\tilde{(k)}} d\rho^{\tilde{(k)}} + \eta \left( \text{TV}(\mu_{h'}) + \text{TV}(\lambda_{h'}) \right).
\end{aligned}
\end{equation}
To the solution of problem \eqref{minpro}, the minimizer $\ell_{h'}^{0} \in \Lambda_{h'}$, obtained by solving problem \eqref{num-minpro}, serves as a numerical approximation. The following theorem concerns the existence of solutions to equation \eqref{num-minpro}, and its proof utilizes Theorem \ref{theory 2} as well as a fact that under the $L^1$ metric, $\Lambda_{h'}$ forms a closed subset of $\Lambda$.
\begin{thm} \label{theory 2}
The minimizing problem \eqref{num-minpro} exists at least one solution $\ell_{h'}^0 \in \Lambda_{h'}$.
\end{thm}
In each subdomain $\Omega^{\tilde{(k)}}$, let us postulate that
\begin{equation} \label{discrete Lame coefficients}
\mu_{h'}(\rho, \phi) = \mathbbm{1}_{(-\infty, 0]}(\rho) \sum_{j = 1}^{m'} \mu_j \mathbbm{1}_{[\varphi_{j - 1},\varphi_j)}(\phi), \quad \lambda_{h'}(\rho, \phi) = \mathbbm{1}_{(-\infty, 0]}(\rho) \sum_{j = 1}^{m'} \lambda_j \mathbbm{1}_{[\varphi_{j - 1}, \varphi_j)}(\phi).
\end{equation}
Then we have
$$
\text{TV}(\mu_{h'}) = \sum_{j=1}^{m'} \tilde{r}(\varphi_j) |\mu_{j+1}-\mu_j|, \quad \text{TV}(\lambda_{h'}) = \sum_{j=1}^{m'} \tilde{r}(\varphi_j) |\lambda_{j+1}-\lambda_j|,
$$
where $\mu_{{m'} + 1} = \mu_1, \lambda_{{m'} + 1} = \lambda_1$. Calculating the gradient of $J$ is crucial in the Adam algorithm, so we will calculate the gradients of $\text{TV}(\mu_{h'})$ and $\text{TV}(\lambda_{h'})$. At the origin, the absolute value $| \cdot |$ is non-differentiability, so we will approximate $\text{TV}(\mu_{h'})$ and $\text{TV}(\lambda_{h'})$:

$$
\text{TV}_{\nu}(\mu_{h'}) = \sum_{j=1}^{m'} \tilde{r}(\varphi_j) \sqrt{|\mu_{j+1}-\mu_j|^2 + \nu^2}, \quad \text{TV}_{\nu}(\lambda_{h'}) = \sum_{j=1}^{m'} \tilde{r}(\varphi_j) \sqrt{|\lambda_{j+1}-\lambda_j|^2 + \nu^2},
$$
where $ \nu \in(0,1]$ is a small regularization parameter. Rather than minimizing $J$ directly, The optimal solution $\ell_{h'}^{\star} \in \Lambda_{h'}$ is obtained by minimizing the regularized energy functional over all candidate parameters $\ell_h' \in \Lambda_{h'}$:
\begin{equation} \label{smoothed energy functional}
\begin{aligned}
J_{\nu}(\ell_h') & := \frac{1}{2} \sum_{\tilde{k}=1}^{\tilde{K}}\sum_{j = 1}^{m'} \int_{-\infty}^0 \int_{\varphi_{j-1 }^{\tilde{(k)}}}^{\varphi_{j}^{\tilde{(k)}}} \left( 2 \mu_{h'} |\varepsilon(u[\ell_{h'}]-z)|^2 + \lambda_{h'} |\nabla \cdot (u[\ell_{h'}]-z)|^2 \right) \\
&e^{2\rho^{\tilde{(k)}}} \tilde{r}^{\tilde{(k)}}(\phi)^2 d\phi^{\tilde{(k)}} d\rho^{\tilde{(k)}} + \eta \left( \text{TV}(\mu_{h'}) + \text{TV}(\lambda_{h'}) \right).
\end{aligned}
\end{equation}
Next, setting $j = 1, \ldots, m'$, we calculate the gradients of $\text{TV}(\mu_{h'})$ and $\text{TV}(\lambda_{h'})$, and $\mu_{0} = \mu_{m'}, \lambda_{0} = \lambda_{m'}$ and $\mu_{{m'} + 1} = \mu_1, \lambda_{{m'} + 1} = \lambda_1$. 
$$
\frac {\partial \text{TV}_{\nu}(\mu_{h'})}{\partial \mu_j} = \tilde{r}(\varphi_{j-1}) \frac{\mu_j-\mu_{j-1}}{\sqrt{|\mu_{j}-\mu_{j - 1}|^2 + \nu^2}} + \tilde{r}(\varphi_{j}) \frac{\mu_{j}-\mu_{j+1}}{\sqrt{|\mu_{j}-\mu_{j + 1}|^2 + \nu^2}},
$$
$$
\frac {\partial \text{TV}_{\nu}(\lambda_{h'})}{\partial \lambda_j} = \tilde{r}(\varphi_{j-1}) \frac{\lambda_j - \lambda_{j-1}}{\sqrt{|\lambda_{j} - \lambda_{j - 1}|^2 + \nu^2}} + \tilde{r}(\varphi_{j}) \frac{\lambda_{j} - \lambda_{j+1}}{\sqrt{|\lambda_{j}-\lambda_{j + 1}|^2 + \nu^2}},
$$

The gradients with respect to the nodal values $\mu_j$ and $\lambda_j$ are computed by applying the chain rule:
$$
\begin{aligned}
\frac {\partial J_0}{\partial \mu_j} & = \int_{-\infty}^0 \int_{0}^{2\pi} \frac{\partial J_0}{\partial \mu_{h'}} \frac{\partial \mu_{h'}}{\partial \mu_j} e^{2\rho^{\tilde{(k)}}} \tilde{r}^{\tilde{(k)}}(\phi)^2 d\phi^{\tilde{(k)}} d\rho^{\tilde{(k)}} = \int_{-\infty}^0 \int_{\varphi_{j-1}}^{\varphi_j} \frac{\partial J_0}{\partial \mu_{h'}} e^{2\rho^{\tilde{(k)}}} \tilde{r}^{\tilde{(k)}}(\phi)^2 d\phi^{\tilde{(k)}} d\rho^{\tilde{(k)}}, \\
\frac {\partial J_0}{\partial \lambda_j} & = \int_{-\infty}^0 \int_{0}^{2\pi} \frac{\partial J_0}{\partial \lambda_{h'}} \frac{\partial \lambda_{h'}}{\partial \lambda_j} e^{2\rho^{\tilde{(k)}}} \tilde{r}^{\tilde{(k)}}(\phi)^2 d\phi^{\tilde{(k)}} d\rho^{\tilde{(k)}} = \int_{-\infty}^0 \int_{\varphi_{j-1}}^{\varphi_j} \frac{\partial J_0}{\partial \lambda_{h'}} e^{2\rho^{\tilde{(k)}}} \tilde{r}^{\tilde{(k)}}(\phi)^2 d\phi^{\tilde{(k)}} d\rho^{\tilde{(k)}},
\end{aligned}
$$
and we know from \cite{part1} that
$$
\frac{\partial J_0}{\partial \mu_{h'}} = -\varepsilon(u[\ell_{h'}]+z)\cdot\varepsilon(u[\ell_{h'}]-z), \quad \frac{\partial J_0}{\partial \lambda_{h'}} = -\frac{1}{2} \left( \nabla \cdot (u[\ell_{h'}]+z) \right) \cdot \left( \nabla \cdot (u[\ell_{h'}]-z) \right).
$$
Hence for $j=1, \ldots, m'$
\begin{equation} \label{gradient of J}
\frac {\partial J_{\nu}}{\partial \mu_j}=\frac {\partial J_0}{\partial \mu_j} + \eta\frac {\partial \text{TV}_{\nu}(\mu_{h'})}{\partial \mu_j}, \quad \frac {\partial J_{\nu}}{\partial \lambda_j}=\frac {\partial J_0}{\partial \lambda_j} + \eta\frac {\partial \text{TV}_{\nu}(\lambda_{h'})}{\partial \lambda_j}.
\end{equation}
After a series of discussions, in order to obtain the gradient and value of $J_{\nu}$, it is necessary to evaluate $u[\ell_h]$. We pay special attention to two issues: First, since $\ell_h \in \Lambda_h$, an interface problem must be considered in the formulation of the forward problem. Therefore, the solution $u[\ell_h]$ exhibits stress characteristics at the intersection points of the interface. Second, the general domains are very complicated with multiple singularities. The direct-line method serves as the forward solver to address both challenges.
Denoting 
$\mathcal{G}^{\tilde{(k)}} = \left( \frac {\partial J_{\nu}}{\partial \mu_1^{\tilde{(k)}}}, \cdots, \frac {\partial J_{\nu}}{\partial \mu_{m'}^{\tilde{(k)}}}, \frac {\partial J_{\nu}}{\partial \lambda_1^{\tilde{(k)}}}, \cdots,\frac {\partial J_{\nu}}{\partial \lambda_{m'}^{\tilde{(k)}}} \right)^T$, 
$\mathscr{L}^{\tilde{(k)}} = (\mu_1^{\tilde{(k)}}, \cdots,\mu_{m'}^{\tilde{(k)}},\lambda_1^{\tilde{(k)}}, \cdots,\lambda_{m'}^{\tilde{(k)}})^T$ ,
and $\mathscr{L}^{\tilde{(k)}}_{t}$ denote the value of $\mathscr{L}^{\tilde{(k)}}$ at the $t$-th iteration, and $\mathcal{G}^{\tilde{(k)}}_t$ denote the value of $\mathcal{G}^{\tilde{(k)}}$ at $\ell_{h'}^t$. And $\ell_{h'}^t = (\mu_{h'}^t, \lambda_{h'}^t)$ represent the discrete Lam\'{e} coefficients defined in \eqref{discrete Lame coefficients} corresponding to $\mathscr{L}^{\tilde{(k)}}_t$. Algorithm \ref{alg: Adam} outlines the integrated computational framework, which incorporates both the Adam optimizer and the direct-line method approach to address the inverse elasticity problem.

\begin{algorithm}[H]
\caption{Implementation of the inverse elasticity problem}
\label{alg: Adam}
\begin{algorithmic}[1]
\State \textbf{Input}: the measurement data $z$, the starting value $\mathscr{L}_0^{\tilde{(k)}}$, the parameters of the regularized functional $J_{\nu}$ (namely $\nu$ and $\eta$), and the algorithmic parameters: exponential decay rates $\hat{\beta}_1, \hat{\beta}_2$, the tolerance threshold $tol$, the learning rate $\tau_t$ at each step, and a small constant $\hat{\epsilon}$.
\State Initialize $t = 0$ and $\mathcal{M}^{\tilde{(k)}}_0 = \mathcal{V}^{\tilde{(k)}}_0 = (0, \ldots, 0)^T \in {R}^{2m}$.
\State Compute the displacement field $u[\ell_{h'}^0]$ by employing the direct-line method, and subsequently evaluate the regularized energy functional $J_{\nu}^0$ using expression \eqref{smoothed energy functional}.
\Do    
\State for $\tilde{k}$= 1,...,$\tilde{K}$    
\State \quad Compute $\mathcal{G}^{\tilde{(k)}}_t$ via \eqref{gradient of J}.    
\State \quad$\mathcal{M}^{\tilde{(k)}}_{t+1} = (1 - \hat{\beta}_1) \mathcal{G}^{\tilde{(k)}}_t + \hat{\beta}_1 \mathcal{M}^{\tilde{(k)}}_t$.    
\State \quad$\mathcal{V}^{\tilde{(k)}}_{t+1} = (1 - \hat{\beta}_2) \mathcal{G}^{\tilde{(k)}}_t \odot \mathcal{G}^{\tilde{(k)}}_t + \hat{\beta}_2 \mathcal{V}^{\tilde{(k)}}_t$, where $\mathcal{G}^{\tilde{(k)}}_t \odot \mathcal{G}^{\tilde{(k)}}_t$ returns the element wise square.    
\State \quad$\tilde{\mathcal{M}}^{\tilde{(k)}}_{t+1} = \mathcal{M}^{\tilde{(k)}}_{t+1} / (1 - (\hat{\beta}_1)^{t + 1}).$    
\State \quad$\tilde{\mathcal{V}}^{\tilde{(k)}}_{t+1} = \mathcal{V}^{\tilde{(k)}}_{t+1} / (1 - (\hat{\beta}_2)^{t + 1}).$    
\State \quad$\mathscr{L}^{\tilde{(k)}}_{t+1} = \mathscr{L}^{\tilde{(k)}}_t - \tau_t \tilde{\mathcal{M}}^{\tilde{(k)}}_{t+1} / (\sqrt{\tilde{\mathcal{V}}^{\tilde{(k)}}_{t+1}} + \hat{\epsilon})$.    
\State Set $t := t + 1$.    
\State Compute the displacement field $u[\ell_{h'}^t]$ by employing the direct-line method, and subsequently evaluate the regularized energy functional $J_{\nu}^0$ using expression \eqref{smoothed energy functional}.
\doWhile{$\frac{| J_{\nu}^t - J_{\nu}^{t - 1} |}{| J_{\nu}^{t - 1} |} > tol.$}
\State \textbf{Output}: $\mathscr{L}_{\tilde{n}}^{\tilde{(k)}}$, $\tilde{n} > 0$ is a certain positive integer.
\end{algorithmic}
\end{algorithm}
The regularized minimization problem \eqref{num-minpro} presents several computational challenges that guide our choice of optimizer. First, the energy functional $J_{\nu}(\ell_h')$ is defined over a high-dimensional parameter space $\Lambda_{h'}$, with the number of parameters scaling with the angular discretization $m'$. Second, the total variation regularization term, while essential for promoting piecewise-constant solutions, introduces a non-smoothness that must be handled. Although the smoothing parameter $\nu$ mitigates this, the objective function can still exhibit complex, non-convex landscapes and varying curvature across different parameters. 
To address these challenges, we employ the Adam optimization algorithm \cite{Kingma2015}. Adam is particularly well-suited for this problem.

\section{Numerical examples} \label{sec:example}

\subsection{Numerical experiments for forward problems} \label{subsec: numerical examples of forward problems}
Numerical solutions and eigenvalues are computed by semi-discretization of linear elements, while reference values are obtained by discretization of quadratic elements on refined meshes in Examples \ref{ex1}-\ref{ex3}. We define the $L^2$ norm $|| \cdot ||_2$ and energy norm $|| \cdot ||_{\star}$ by
$$
|| v ||_2 = \left( \iint_{\Omega} \left( |v_1|^2+|v_2|^2 \right) dxdy \right)^{1/2}, \quad || v ||_{\star} = \left(  \iint_{\Omega} \left( 2\mu|\varepsilon(v)|^2 + \lambda|\nabla \cdot v|^2 \right) dxdy \right)^{1/2},
$$

We can conclude that for any $v \in H^1(\Omega) \times H^1(\Omega)$, the inequality $|| v ||_{\star} \le || v ||_{*}$ holds, where $|| \cdot ||_{*}$ has been defined in Theorem \ref{thm-error}. We will judge the feasibility of our model by the energy relative error $\frac{||u-u^h||_{\star}}{||u||_{\star}}$ and $L^2$ relative error $\frac{||u-u^h||_2}{||u||_2}$, where $u^h = (u_1^h, u_2^h)^T$ is the numerical solution measured using the direct-line method and $u = (u_1, u_2)^T$ is approximated as the reference/exact solution.

\begin{example}\label{ex1}
\end{example}
Let the boundary functions be defined as:
$$
\tilde{r}^{(1)}(\phi) =
\begin{cases}
\dfrac{1}{\sqrt{\sin^4(\phi) + \cos^4(\phi)}}, & \quad  \pi/4 < \phi \le 7\pi/4,\\
-\dfrac{1}{\sin(\phi-\frac{\pi}{2})}, & \quad -\pi/4 < \phi \le \pi/4,
\end{cases}
$$
$$
\tilde{r}^{(2)}(\phi) =
\begin{cases}
\dfrac{1}{\sqrt{\sin^4(\phi) + \cos^4(\phi)}}, & \quad  -3\pi/4 < \phi \le 3\pi/4,\\
-\dfrac{1}{\sin(\phi+\frac{\pi}{2})}, & \quad 3\pi/4 < \phi \le 5\pi/4,
\end{cases}
$$
with $\Omega^{(\tilde{k})} = \{ (r, \phi)| 0 \le r < \tilde{r}^{(\tilde{k})}(\phi), 0 \le \phi < 2\pi, \tilde{k}=1,2 \}$, $\Omega_{i}^{(\tilde{k})} = \{ (r, \phi)| 0 < r < \tilde{r}^{(\tilde{k})}(\phi),  \theta_i^{(\tilde{k})} < \phi < \theta_{i + 1}^{(\tilde{k})}, i = 1, 2\}$, $\theta_1^{(1)} = -\pi, \theta_2^{(1)} = -\pi/4, \theta_3^{(1)} = \pi$, and $\theta_1^{(2)} = 0, \theta_2^{(2)} = 3\pi/4, \theta_3^{(2)} = 2\pi$. 

The lam\'e coefficients here are $(\mu_1^{(1)}, \lambda_1^{(1)}) = (0.1, 1)$, $(\mu_2^{(1)}, \lambda_2^{(1)}) = (0.6, 2)$, $(\mu_1^{(2)}, \lambda_1^{(2)}) = (0.1, 1)$, $(\mu_2^{(2)}, \lambda_2^{(2)}) = (0.6, 2)$; And the body force $p=(1,1)^T$ and the Dirichlet boundary condition $f=(y^2,-x^2)^T$. see Figure \ref{fig:domainExp1-1}.
\begin{figure}[H]
\centering
\includegraphics[scale=0.6]{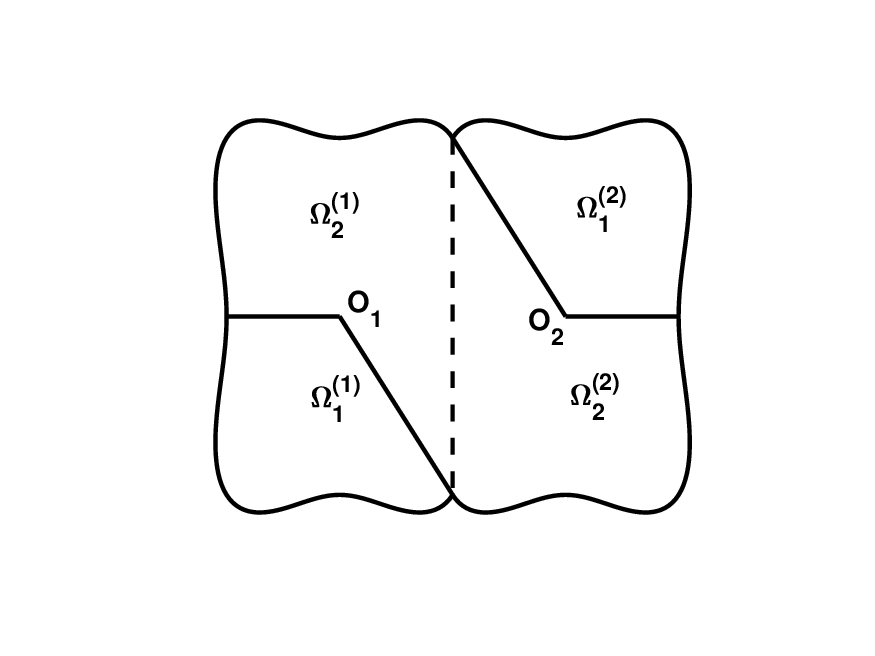}
\caption{Domain $\Omega$ in Example \ref{ex1}}
\label{fig:domainExp1-1}
\end{figure}
We show $\frac{\partial u}{\partial r}(r, \phi_0)$ with $\phi_0 = \frac{3\pi}{4}$ in Figure \ref{singularity of ex1-1}. It can be observed that as $r \rightarrow 0$, the derivatives $\frac{\partial u_1^{(1)}}{\partial r}$ and $\frac{\partial u_2^{(1)}}{\partial r}$ tend to $-\infty$ and $\infty$, respectively. This behavior confirms the presence of stress singularities at the singular point.

\begin{figure}[H]
\centering
\begin{subfigure}[b]{0.4\textwidth}
\includegraphics[width=\textwidth]{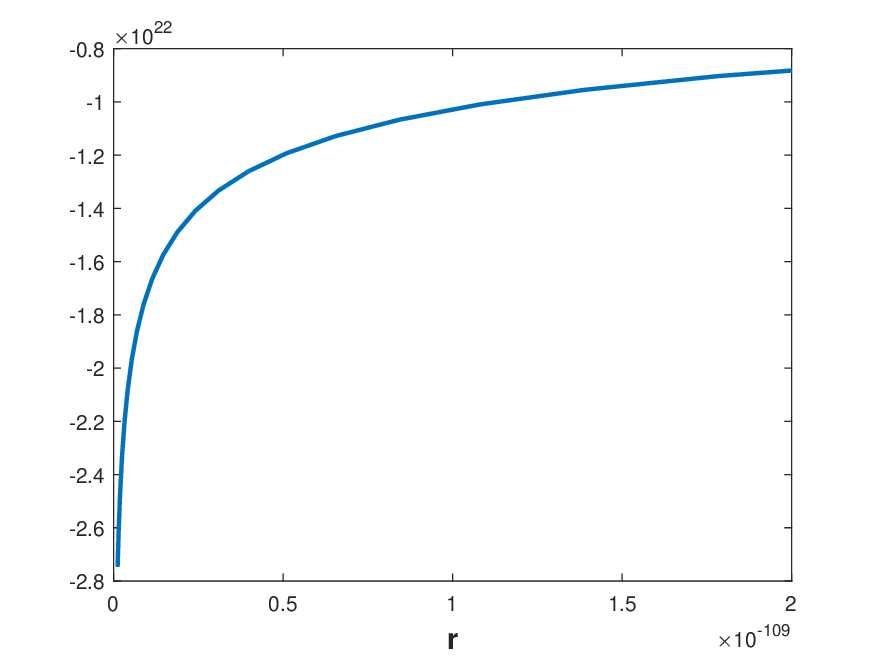}
\caption{$\frac{\partial u_1^{(1)}}{\partial r}$}
\end{subfigure}
\begin{subfigure}[b]{0.4\textwidth}
\includegraphics[width=\textwidth]{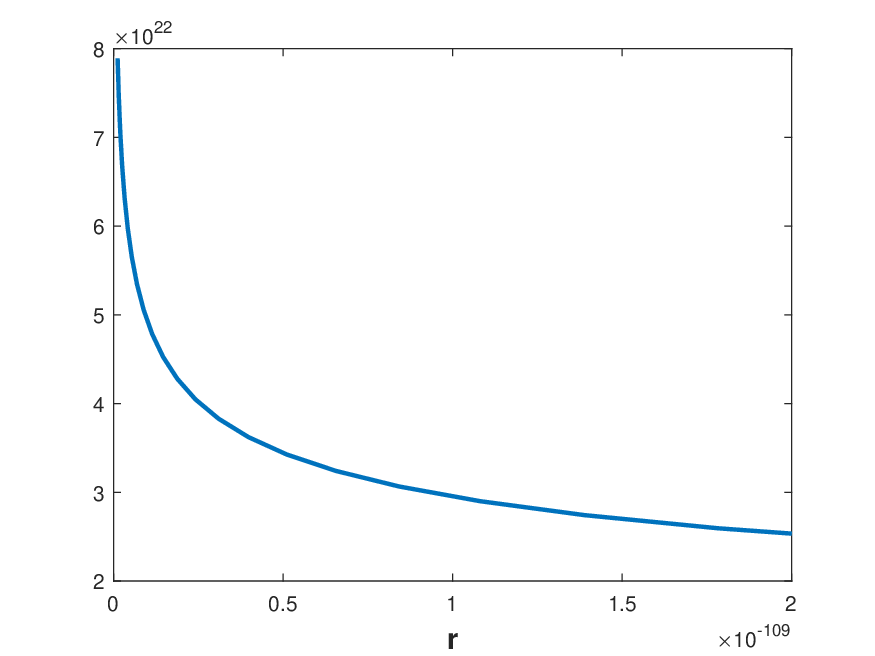}
\caption{$\frac{\partial u_2^{(1)}}{\partial r}$}
\end{subfigure}
\caption{$\frac{\partial u}{\partial r}(r, \phi_0)$ with $\phi_0 = \frac{3\pi}{4}$ in Example \ref{ex1}}
\label{singularity of ex1-1}
\end{figure}
The errors of the first non-zero eigenvalue $\gamma^h_3$ relative to the reference value $\gamma_3=0.7775070808$ are provided in Table \ref{tab:multi-material1-1}. As the mesh refines, exhibiting a second-order convergence rate. This result indicates that the eigenvalues can achieve quadratic convergence towards exact eigenvalues, proving the method's natural ability to capture singularities.
The $L^2$ and energy relative errors are shown in Table \ref{tab:errors-multi-materials1-2}. The $L^2$ error converges at a rate of 2 , while the energy error converges linearly. These results align with Theorem \ref{thm-error}, which guarantees first-order convergence in the energy norm and second-order in the $L^2$ norm for linear elements, even in the presence of stress singularities.
\begin{table}[H]
\centering
\caption{Errors of $\gamma_3^h$, with $\gamma_3 =  0.7775070808$}
\label{tab:multi-material1-1}
\begin{tabular}{ccc}
\toprule
M & $| \gamma_3^h - \gamma_3 |$ & order \\
\midrule
16 
& 1.638e-2 &  --  \\
32 
& 3.835e-3 & 2.094 \\
64 
& 9.181e-4 & 2.063 \\
128 
& 2.268e-4 & 2.018 \\
256 
& 5.654e-5 & 2.005 \\
\bottomrule
\end{tabular}
\end{table}

\begin{table}[H]  
\centering  
\caption{Convergence order of relative error between $L^2$ and energy}  
\label{tab:errors-multi-materials1-2}  
\begin{tabular}{ccccc}
\toprule
M & $\frac{\|u-u^h\|_2}{\|u\|_2}$ & order & $\frac{\|u-u^h\|_{\star}}{\|u\|_{\star}}$ & order \\    
\midrule
16 
& 7.229e-2 & --         & 8.208e-2 & --         \\
32 
& 2.304e-2 & 1.649      & 3.324e-2 & 1.304      \\
64 
& 6.277e-3 & 1.877      & 1.497e-2 & 1.151      \\
128 
& 1.660e-3 & 1.919     & 7.094e-3 & 1.078      \\
256 
& 4.552e-4 & 1.867     & 3.415e-3 & 1.055      \\    
\bottomrule
\end{tabular}
\end{table}

\begin{example}\label{ex2}
\end{example}
Let the boundary functions be defined as:
$$
\tilde{r}^{(1)}(\phi) =
\begin{cases}
\dfrac{1}{\sqrt{\sin^4(\phi) + \cos^4(\phi)}}, & \quad  -\pi/4 < \phi \le 3\pi/4,\\
-\dfrac{1}{\sin(\phi+\frac{\pi}{2})}, & \quad 3\pi/4 < \phi \le 5\pi/4,\\
-\dfrac{1}{\sin(\phi)}, & \quad -3\pi/4 < \phi \le -\pi/4,
\end{cases}
$$
 $$
\tilde{r}^{(2)}(\phi) =
\begin{cases}
\dfrac{1}{\sqrt{\sin^4(\phi) + \cos^4(\phi)}}, & \quad  \pi/4 < \phi \le 5\pi/4,\\
-\dfrac{1}{\sin(\phi)}, & \quad -3\pi/4 < \phi \le -\pi/4,\\
-\dfrac{1}{\sin(\phi-\frac{\pi}{2})}, & \quad -\pi/4 < \phi \le \pi/4,
\end{cases}
$$

$$
\tilde{r}^{(3)}(\phi) =
\begin{cases}
\dfrac{1}{\sqrt{\sin^4(\phi) + \cos^4(\phi)}}, & \quad  3\pi/4 < \phi \le 7\pi/4,\\
-\dfrac{1}{\sin(\phi-\frac{\pi}{2})}, & \quad -\pi/4 < \phi \le \pi/4,\\
-\dfrac{1}{\sin(\phi-\pi)}, & \quad \pi/4 < \phi \le 3\pi/4,
\end{cases}
$$

$$
\tilde{r}^{(4)}(\phi) =
\begin{cases}
\dfrac{1}{\sqrt{\sin^4(\phi) + \cos^4(\phi)}}, & \quad  -3\pi/4 < \phi \le \pi/4,\\
-\dfrac{1}{\sin(\phi-\pi)}, & \quad \pi/4 < \phi \le 3\pi/4,\\
-\dfrac{1}{\sin(\phi+\frac{\pi}{2})}, & \quad 3\pi/4 < \phi \le 5\pi/4,
\end{cases}
$$
with $\Omega_{i}^{(\tilde{k})}= \{ (r, \phi)| 0 < r < \tilde{r}^{(\tilde{k})}(\phi),  \theta_i^{(\tilde{k})} < \phi < \theta_{i + 1}^{(\tilde{k})}, \tilde{k}=1,...,4,i=1,2\}$,  $\theta_1^{(1)} = 0, \theta_2^{(1)} = \pi/2, \theta_3^{(1)} = 2\pi$;  $\theta_1^{(2)} = -\pi, \theta_2^{(2)} = \pi/2, \theta_3^{(2)} = \pi$;  $\theta_1^{(3)} = -\pi, \theta_2^{(3)} = -\pi/2, \theta_3^{(3)} = \pi$; and $\theta_1^{(4)} = 0, \theta_2^{(4)} = 3\pi/2, \theta_3^{(4)} = 2\pi$. 

The lam\'e coefficients here are $(\mu_1^{(1)}, \lambda_1^{(1)}) = (0.6, 2)$, $(\mu_2^{(1)}, \lambda_2^{(1)}) = (0.1, 1)$; $(\mu_1^{(2)}, \lambda_1^{(2)}) = (0.7, 2.1)$, $(\mu_2^{(2)}, \lambda_2^{(2)}) = (0.1, 1)$; $(\mu_1^{(3)}, \lambda_1^{(3)}) = (0.6, 2)$, $(\mu_2^{(3)}, \lambda_2^{(3)}) = (0.1, 1)$; $(\mu_1^{(4)}, \lambda_1^{(4)}) = (0.5, 1.9)$, $(\mu_2^{(4)}, \lambda_2^{(4)}) = (0.1, 1) $; And the body force $p=(1,1)^T$ and the Dirichlet boundary condition $f=(1,1)^T$. see Figure \ref{fig:domainExp2-1}. 
\begin{figure}[H]
\centering
\includegraphics[scale=0.6]{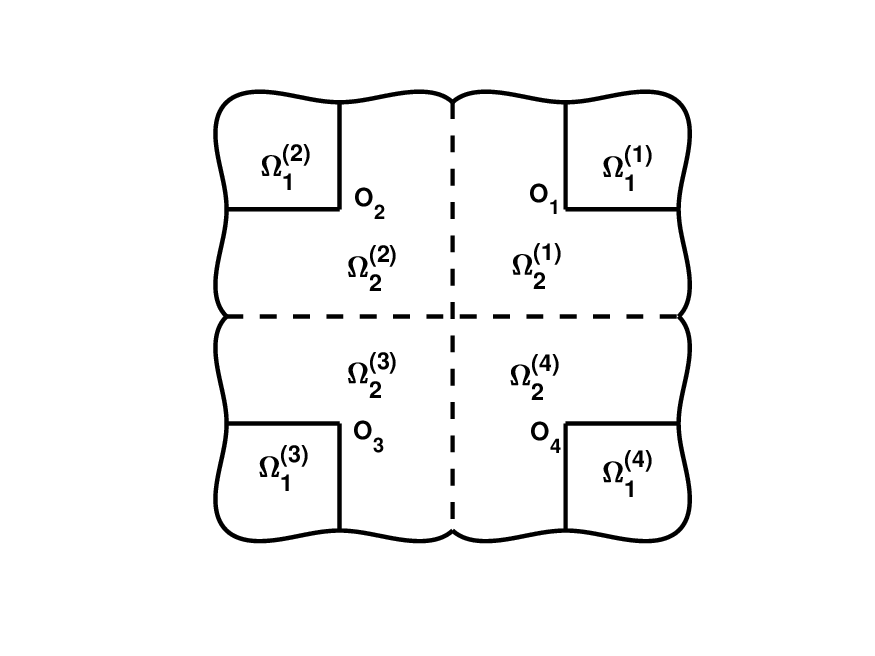}
\caption{Domain $\Omega$ in Example \ref{ex2}}
\label{fig:domainExp2-1}
\end{figure}
We show $\frac{\partial u}{\partial r}(r, \phi_0)$ with $\phi_0 = \frac{\pi}{4}$ in Figure \ref{singularity of ex2-1}. It can be observed that as $r \rightarrow 0$, the derivatives $\frac{\partial u_1^{(1)}}{\partial r}$ and $\frac{\partial u_2^{(1)}}{\partial r}$ tend to $-\infty$ and $\infty$, respectively. This behavior confirms the presence of stress singularities at the singular point.

\begin{figure}[H]
\centering
\begin{subfigure}[b]{0.4\textwidth}
\includegraphics[width=\textwidth]{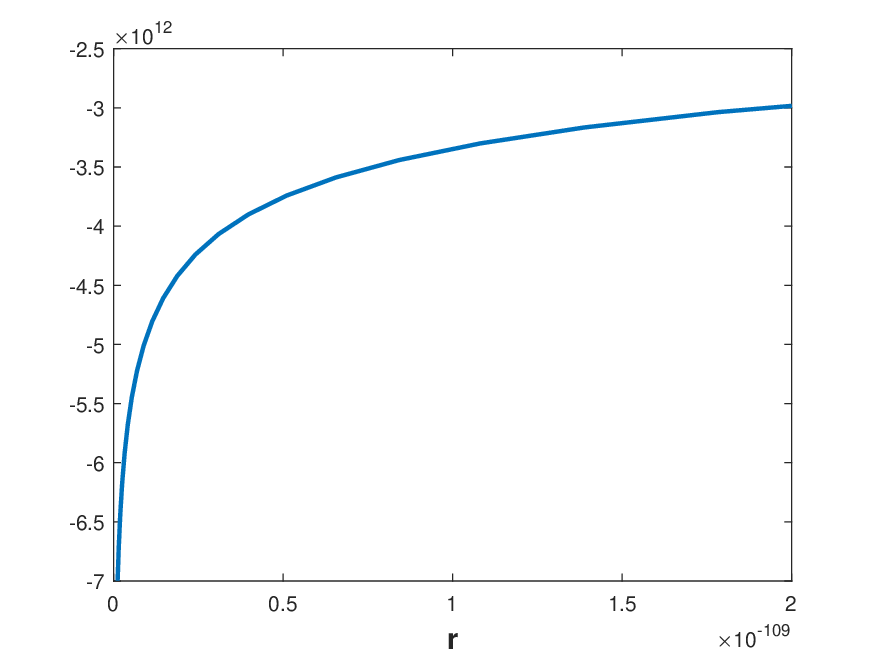}
\caption{$\frac{\partial u_1^{(1)}}{\partial r}$}
\end{subfigure}
\begin{subfigure}[b]{0.4\textwidth}
\includegraphics[width=\textwidth]{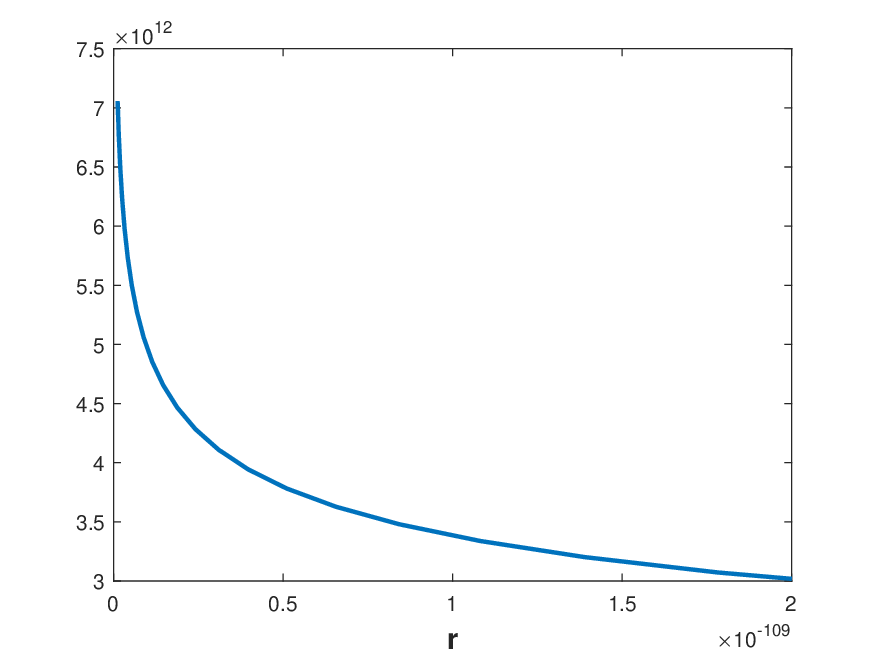}
\caption{$\frac{\partial u_2^{(1)}}{\partial r}$}
\end{subfigure}
\caption{$\frac{\partial u}{\partial r}(r, \phi_0)$ with $\phi_0 = \frac{\pi}{4}$ in Example \ref{ex2}}
\label{singularity of ex2-1}
\end{figure}
The errors of the first non-zero eigenvalue $\gamma^h_3$ relative to the reference value $\gamma_3=0.8333347413$ are provided in Table \ref{tab:multi-material2-1}. As the mesh refines, exhibiting a second-order convergence rate. This result indicates that the eigenvalues can achieve quadratic convergence towards exact eigenvalues, proving the method's natural ability to capture singularities.
The $L^2$ and energy relative errors are shown in Table \ref{tab:errors-multi-materials2-1}. The $L^2$ error converges at a rate of 2 , while the energy error converges linearly. These results align with Theorem \ref{thm-error}, which guarantees first-order convergence in the energy norm and second-order in the $L^2$ norm for linear elements, even in the presence of stress singularities.
\begin{table}[H]
\centering
\caption{Errors of $\gamma_3^h$, with $\gamma_3 =  0.8333347413$}
\label{tab:multi-material2-1}
\begin{tabular}{ccc}
\toprule
{M} & {$| \gamma_3^h - \gamma_3 |$} & {order} \\
\midrule
16  
& 5.594e-2 & --   \\
32  
& 3.173e-2 & 0.818  \\
64  
& 1.035e-2 & 1.616  \\
128 
& 2.811e-3 & 1.881  \\
256 
& 7.186e-4 & 1.968  \\
\bottomrule
\end{tabular}
\end{table}

\begin{table}[H]  
\centering
\caption{Convergence order of relative error between $L^2$ and energy}  
\label{tab:errors-multi-materials2-1}  
\begin{tabular}{ccccc}
\toprule
M 
& $\frac{\|u-u^h\|_2}{\|u\|_2}$ & order & $\frac{\|u-u^h\|_{\star}}{\|u\|_{\star}}$ & order \\    
\midrule
16  
& 3.441e-2 & --          & 3.795e-2 & --          \\
32  
& 1.215e-2 & 1.502       & 1.538e-2 & 1.303       \\
64  
& 3.267e-3 & 1.896       & 5.991e-3 & 1.361       \\
128 
& 8.397e-4 & 1.961       & 2.678e-3 & 1.163       \\
256 
& 2.139e-4 & 1.973       & 1.293e-3 & 1.051       \\    
\bottomrule
\end{tabular}
\end{table}

\begin{example}\label{ex3}
\end{example}
Here, the middle is empty, $ \Omega= \bigcup \Omega_k,k=1,..,5$, $\Omega = \bigcup \Omega^{(\tilde{k})}, \tilde{k}=1,...,4$, and $\Omega^{(\tilde{k})}=\bigcup\Omega^{(\tilde{k})}_i,i=1,2,3$. see Figure \ref{fig:domain33}. 
\begin{figure}[H]
\centering
\begin{subfigure}[b]{0.46\textwidth}  
\includegraphics[width=\textwidth]{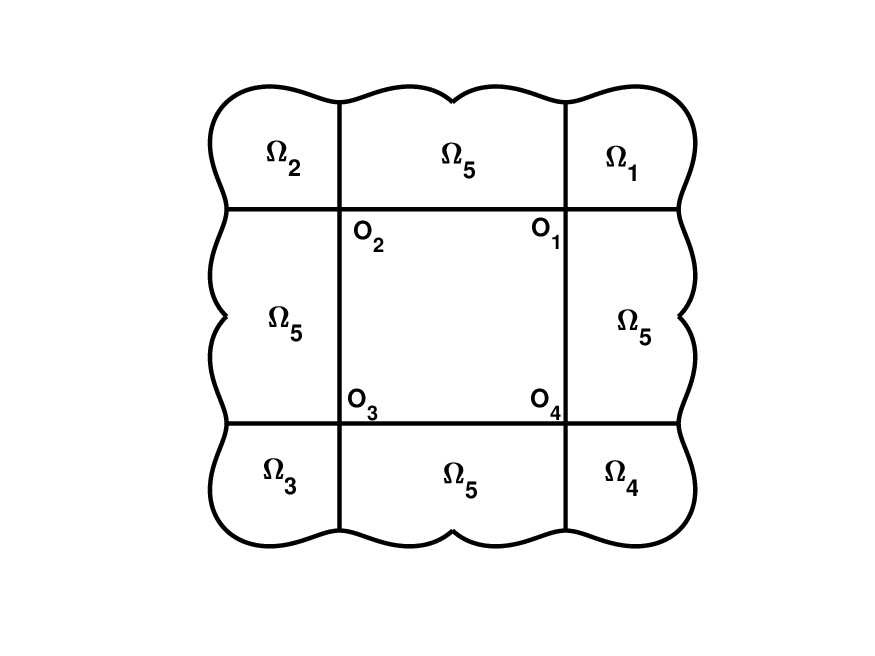}  
\caption{Domain $\Omega$}  
\label{fig:domain3a}
\end{subfigure}
\begin{subfigure}[b]{0.46\textwidth}  
\includegraphics[width=\textwidth]{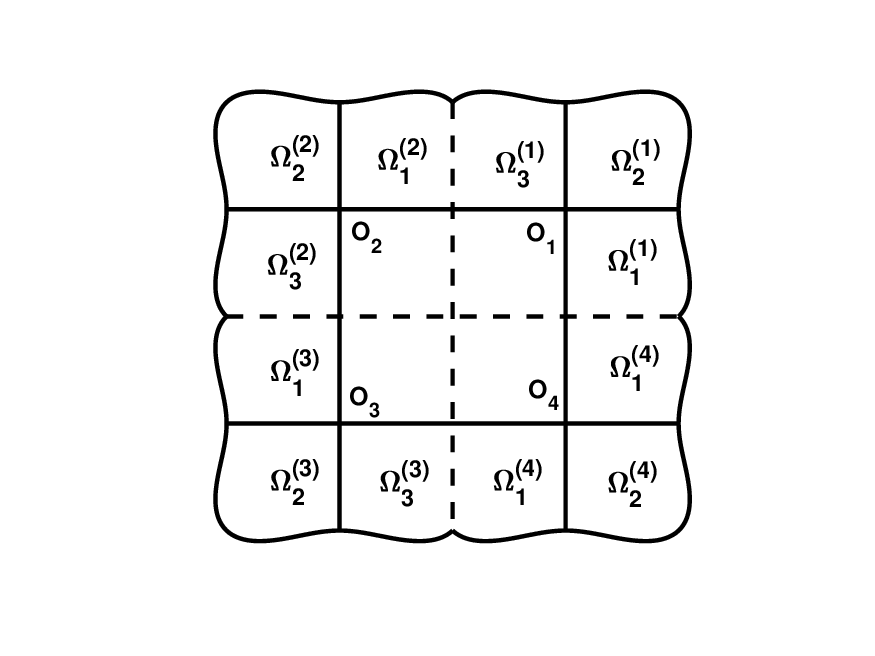}  
\caption{Subdomain $\Omega^{(\tilde{k})}$}  
\label{fig:domain3b}
\end{subfigure}
\caption{Decomposition of domain $\Omega$ into subdomains with a central void and four singular points.}
\label{fig:domain33}
\end{figure}
Let the boundary functions be defined as:
$$
\tilde{r}^{(1)}(\phi) =
\begin{cases}
-\dfrac{1}{\sin(\phi)}, & \quad -\pi/2 \le \phi \le -\pi/4,\\
\sqrt{1 + \sin^2(2\phi)}, & \quad  -\pi/4 < \phi \le 3\pi/4,\\
-\dfrac{1}{\sin(\phi+\frac{\pi}{2})}, & \quad 3\pi/4 < \phi \le \pi,\\
\end{cases}
$$

$$
\tilde{r}^{(2)}(\phi) =
\begin{cases}
-\dfrac{1}{\sin(\phi-\frac{\pi}{2})}, & \quad 0 \le \phi \le \pi/4,\\
\sqrt{1 + \sin^2(2\phi)}, & \quad  \pi/4 < \phi \le 5\pi/4,\\
-\dfrac{1}{\sin(\phi)}, & \quad -3\pi/4 < \phi \le -\pi/2,\\
\end{cases}
$$

$$
\tilde{r}^{(3)}(\phi) =
\begin{cases}
-\dfrac{1}{\sin(\phi-\pi)}, & \quad \pi/2 \le \phi \le 3\pi/4,\\
\sqrt{1 + \sin^2(2\phi)}, & \quad  3\pi/4 < \phi \le 7\pi/4,\\
-\dfrac{1}{\sin(\phi-\frac{\pi}{2})}, & \quad -\pi/4 < \phi \le 0,\\
\end{cases}
$$

$$
\tilde{r}^{(4)}(\phi) =
\begin{cases}
-\dfrac{1}{\sin(\phi+\frac{\pi}{2})}, & \quad -\pi \le \phi \le -3\pi/4,\\
\sqrt{1 + \sin^2(2\phi)}, & \quad  -3\pi/4 < \phi \le \pi/4,\\
-\dfrac{1}{\sin(\phi-\pi)}, & \quad \pi/4 < \phi \le \pi/2,\\
\end{cases}
$$
with $\Omega_{i}^{( \tilde{k})}= \{ (r, \phi)| 0 < r < \tilde{r}^{(\tilde{k})}(\phi),  \theta_i^{(\tilde{k})} < \phi < \theta_{i + 1}^{(\tilde{k})}, \tilde{k} =1,...,4,i=1,2,3\}$,  $\theta_1^{(1)} = -\pi/2, \theta_2^{(1)} =0, \theta_3^{(1)}=\pi/2,\theta_4^{(1)} =\pi$;  $\theta_1^{(2)} = 0, \theta_2^{(2)} =\pi/2, \theta_3^{(2)}=\pi,\theta_4^{(2)} =3\pi/2$; $\theta_1^{(3)} = \pi/2, \theta_2^{(3)} =\pi, \theta_3^{(3)}=3\pi/2,\theta_4^{(3)} =2\pi$; and $\theta_1^{(4)} = -\pi, \theta_2^{(4)} =-\pi/2, \theta_3^{(4)}=0,\theta_4^{(4)} =\pi/2$. 

The lam\'e coefficients here are $(\mu_1^{(1)}, \lambda_1^{(1)}) = (0.4, 1)$, $(\mu_2^{(1)}, \lambda_2^{(1)}) = (0.6, 1.5)$, $(\mu_3^{(1)}, \lambda_3^{(1)}) = (0.4, 1) $;  $(\mu_1^{(2)}, \lambda_1^{(2)}) = (0.4, 1)$, $(\mu_2^{(2)}, \lambda_2^{(2)}) = (0.6, 1.5)$, $(\mu_3^{(2)}, \lambda_3^{(2)}) = (0.4, 1) $;  $(\mu_1^{(3)}, \lambda_1^{(3)}) = (0.4, 1)$, $(\mu_2^{(2)}, \lambda_2^{(3)}) = (0.6, 1.5)$, $(\mu_3^{(3)}, \lambda_3^{(3)}) = (0.4, 1) $;  $(\mu_1^{(4)}, \lambda_1^{(4)}) = (0.4, 1)$, $(\mu_2^{(4)}, \lambda_2^{(4)}) = (0.6, 1.5)$, $(\mu_3^{(4)}, \lambda_3^{(4)}) = (0.4, 1) $; And the body force $p=(\frac{x}{\sqrt{x^2+y^2}},\frac{y}{\sqrt{x^2+y^2}})^T$ and the Dirichlet boundary condition $f=(1,1)^T$.

We show $\frac{\partial u}{\partial r}(r, \phi_0)$ with $\phi_0 = \frac{3\pi}{4}$ in Figure \ref{singularity of ex3}. It can be observed that as $r \rightarrow 0$, the derivatives $\frac{\partial u_1^{(2)}}{\partial r}$ and $\frac{\partial u_2^{(2)}}{\partial r}$ tend to $-\infty$ and $\infty$, respectively. This behavior confirms the presence of stress singularities at the singular point.

\begin{figure}[H]
\centering
\begin{subfigure}[b]{0.4\textwidth}
\includegraphics[width=\textwidth]{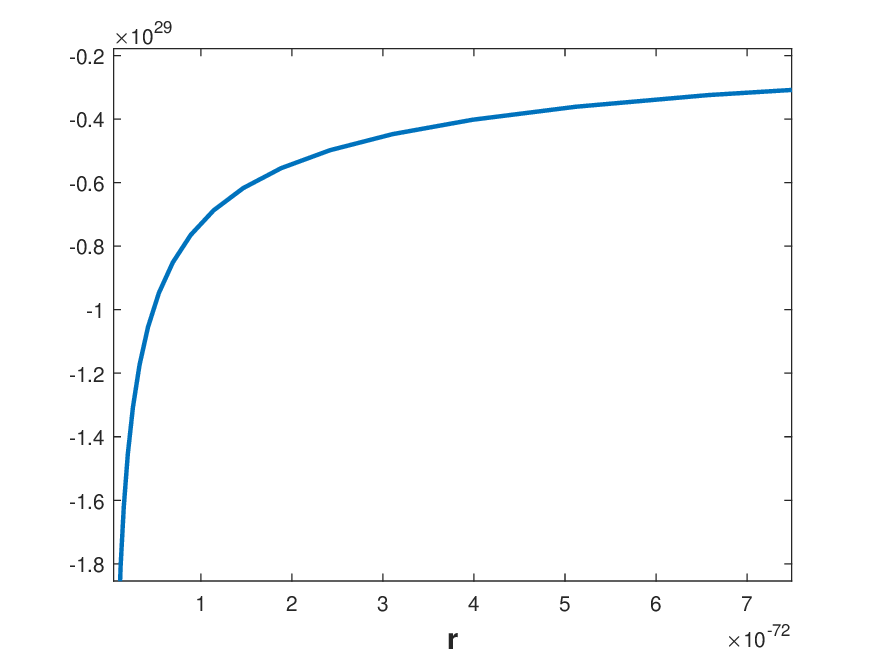}
\caption{$\frac{\partial u_1^{(2)}}{\partial r}$}
\end{subfigure}
\begin{subfigure}[b]{0.4\textwidth}
\includegraphics[width=\textwidth]{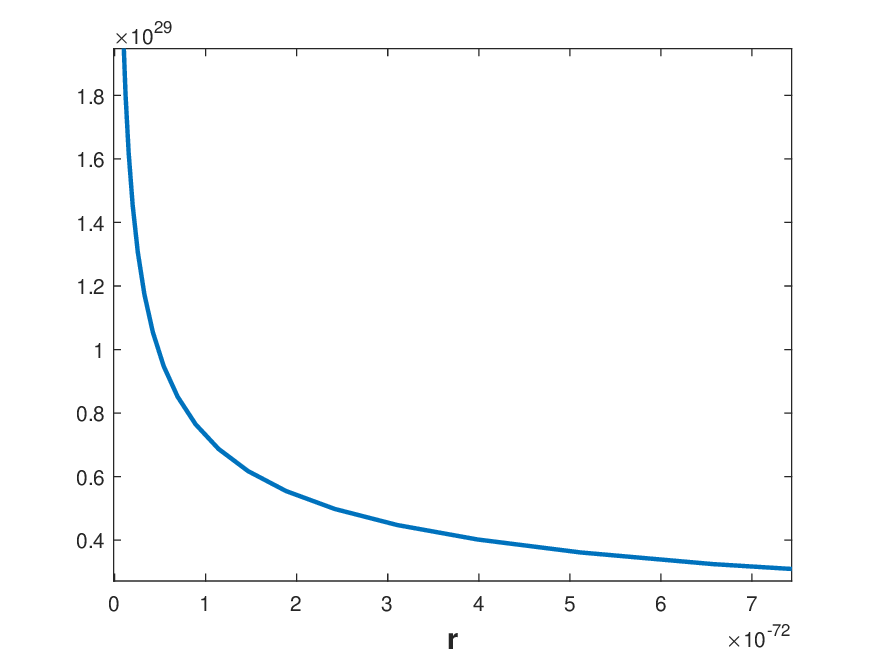}
\caption{$\frac{\partial u_2^{(2)}}{\partial r}$}
\end{subfigure}
\caption{$\frac{\partial u}{\partial r}(r, \phi_0)$ with $\phi_0 = \frac{3\pi}{4}$ in Example \ref{ex1}}
\label{singularity of ex3}
\end{figure}
The errors of the first non-zero eigenvalue $\gamma^h_3$ relative to the reference value $\gamma_3=0.5713189037$ are provided in Table \ref{tab:multi-material}. As the mesh refines, exhibiting a second-order convergence rate. This result indicates that the eigenvalues can achieve quadratic convergence towards exact eigenvalues, proving the method's natural ability to capture singularities.
The $L^2$ and energy relative errors are shown in Table \ref{tab:errors-multi-materials} . The $L^2$ error converges at a rate of 2 , while the energy error converges linearly. These results align with Theorem \ref{thm-error}, which guarantees first-order convergence in the energy norm and second-order in the $L^2$ norm for linear elements, even in the presence of stress singularities.
\begin{table}[H]
\centering
\caption{Errors of $\gamma_3^h$, with $\gamma_3 =  0.5713189037$}
\label{tab:multi-material}
\begin{tabular}{ccc}
\toprule
M 
& $| \gamma_3^h - \gamma_3 |$ & order \\
\midrule
16 
& 4.846e-3 &  \\
32 
& 1.308e-3 & 1.890 \\
64 
& 3.305e-4 & 1.985 \\
128 
& 8.287e-5 & 1.997 \\
256 
& 2.073e-5 & 1.999 \\
\bottomrule
\end{tabular}
\end{table}

\begin{table}[H]  
\centering
\caption{Convergence order of relative error between $L^2$ and energy}  
\label{tab:errors-multi-materials}  
\begin{tabular}{ccccc}
\toprule
{M} 
& {$\frac{\|u-u^h\|_2}{\|u\|_2}$} & {order} & {$\frac{\|u-u^h\|_{\star}}{\|u\|_{\star}}$} & {order} \\    
\midrule
16  
& 1.339e-3 & --       & 3.501e-3 & --       \\
32  
& 5.161e-4 & 1.375    & 1.899e-3 & 0.882    \\
64  
& 1.409e-4 & 1.873    & 9.838e-4 & 0.949    \\
128 
& 3.598e-5 & 1.969    & 4.972e-4 & 0.984    \\
256 
& 8.835e-6 & 2.026    & 2.491e-4 & 0.997    \\    
\bottomrule
\end{tabular}
\end{table}

Numerical experiments (examples \ref{ex1}-\ref{ex3}) show that the eigenvalues using direct-line method converge rapidly to the true eigenvalue of elliptic operator, confirming its inherent ability of this method to solve for singularities. Notably, even under the condition that there are multiple stress singularities in general regions, optimal convergence rates are achieved for linear/quadratic elements. These results collectively validate the effectiveness and precision of the proposed approach for modeling linear elasticity in composite materials within general domains.
\subsection{Numerical experiments for inverse elasticity problems}
Next, our numerical experiment use the Adam Algorithm \ref{alg: Adam} to minimize \eqref{smoothed energy functional} over $\ell_{h'} \in \Lambda_{h'}$, with some parameter settings as follows: $tol = 5 \times 10^{-6},\hat{\beta}_1 = 0.9, \hat{\beta}_2 = 0.999, \eta = 10^{-7}, \nu = 10^{-7},  \hat{\epsilon} = 10^{-7}$ and a decaying $\tau_k$. And set the grid size $h'=\pi/64$ to discretize the solution space, that is 
$$\mu_{h'} = \mathbbm{1}_{(-\infty, 0]}(\rho) \sum_{j=1}^{128} \mu_j \mathbbm{1}_{[(j-1)h', jh')}(\phi),
\quad \lambda_{h'} = \mathbbm{1}_{(-\infty, 0]}(\rho) \sum_{j=1}^{128} \lambda_j \mathbbm{1}_{[(j-1)h', jh')}(\phi).$$
Here, we use $\mu_j^k$ and $\lambda_j^k$ to represent the values of $\mu_j$ and $\lambda_j$ in step $k$ respectively. At each iteration, each iteration is solving a forward problem (mesh parameter $h'=\pi/64$). For validation purposes, reference solutions $u[\ell^{\star}]$ are computed using quadratic elements on refined meshes via the same direct-line method, providing noise-free benchmark results. The measurement data are collected at discrete points within the domain $\Omega$, denoted by the set $\Xi = \{ (x_j, y_j) \}_{j = 1}^{M_1} \subset \Omega$. To quantify the precision of the reconstructed Lam\'{e} coefficients, the $L^1$ error is used, where the $L^1$ norm is represented by $|| \cdot ||_1$.

\begin{example}\label{example41}
\end{example}
Building on the analysis of the forward problem in Example \ref{ex1}, we now investigate the corresponding inverse elasticity problem to estimate the Lam\'{e} coefficients from displacement measurements. 

The noisy measurement is defined as $z|_{\Xi} =\delta \xi_U || u[\ell^{\star}] ||_{\infty} +u[\ell^{\star}]|_{\Xi}$, where $|| \cdot ||_{\infty}$ denotes in $L^{\infty}$ norm, $\delta=0.0001$ represents the noise level and $\xi_U$ is a vector of independent identically distributed uniform random variables on [-1,1]. The initial values are set as $\mu_1^0 = \cdots = \mu_{128}^0 = 0.35, \lambda_1^0 = \cdots = \lambda_{128}^0 = 1.5$. The iteration process continues for 842 steps until the termination criterion is met. Figure \ref{fig:J1} reflects the change in the value of $J$, which converges gradually with the iteration step number $k$. The numerically reconstructed Lam\'{e} coefficients $\mu_{h'}^{842}$ and $\lambda_{h'}^{842}$ along the $\phi$-direction is represented by the blue broken line in Figure \ref{one-lam1-mu41}, and the true values are indicated by the red dotted lines. In the 842nd iteration, the error results are $\|\mu_{h'}^{842} - \mu^{\star}\|_{1}/\|\mu^{\star}\|_{1} = \text{7.873e-3}$ and $\|\lambda_{h'}^{842} - \lambda^{\star}\|_{1}/\|\lambda^{\star}\|_{1} = \text{1.111e-2}$.
\begin{figure}[H]
\centering
\includegraphics[scale=0.4]{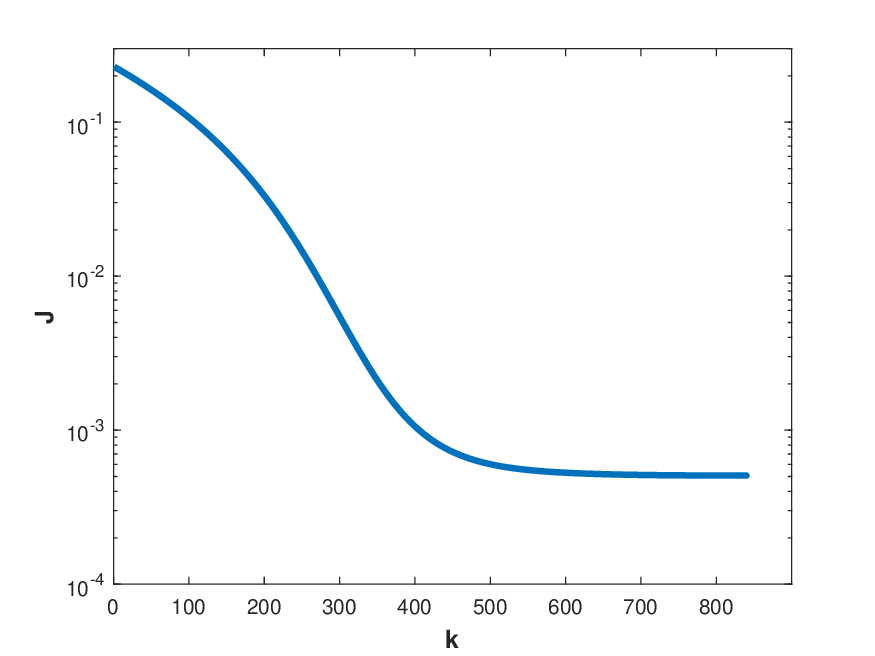}
\caption{Convergence of the objective function $J$ in Example \ref{example41}} 
\label{fig:J1}
\end{figure}
\begin{figure}[H]
\centering
\begin{subfigure}[b]{0.4\textwidth}
\includegraphics[width=\textwidth]{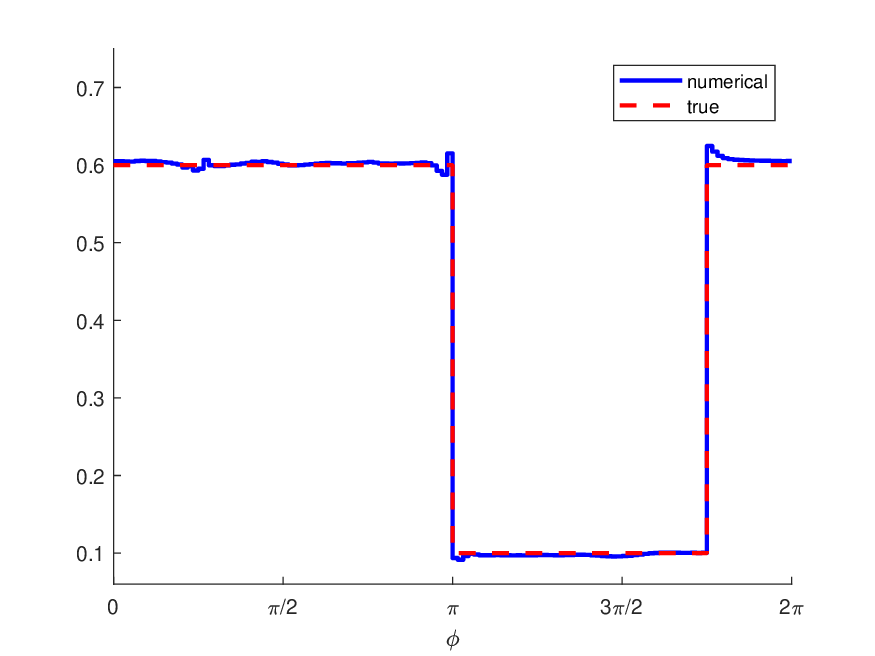}
\caption{$\mu$}
\end{subfigure}
\begin{subfigure}[b]{0.4\textwidth}
\includegraphics[width=\textwidth]{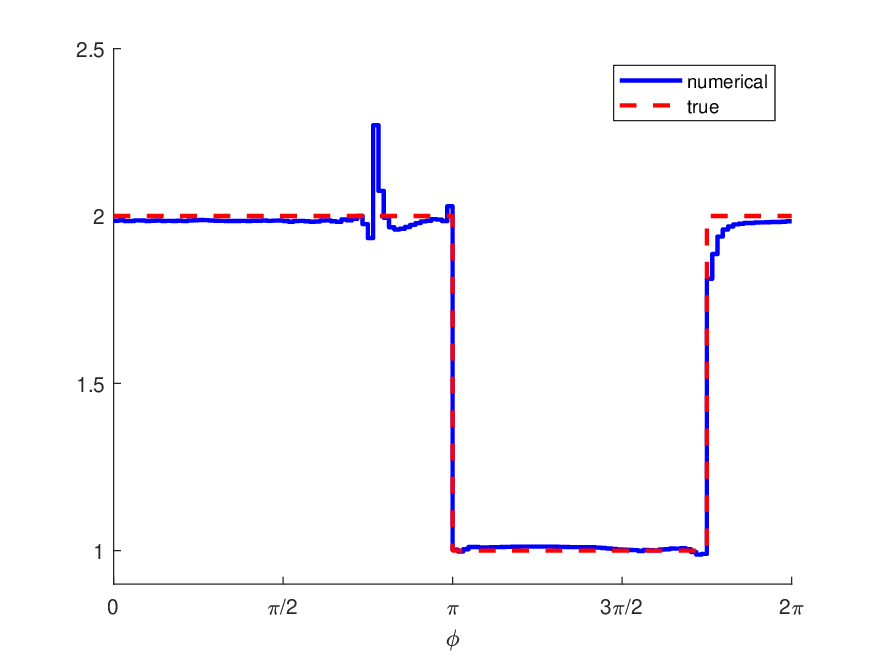}
\caption{$\lambda$}
\end{subfigure}
\caption{Reconstructed Lam\'e coefficients and true Lam\'e coefficients}
\label{one-lam1-mu41}
\end{figure}
\begin{example}\label{example2}
\end{example}
Building upon the forward problem analysis in Example \ref{ex2}, we now investigate the corresponding inverse elasticity problem to estimate the Lam\'{e} coefficients from displacement measurements.

The noisy measurement is defined as $z|_{\Xi} =\delta \xi_U || u[\ell^{\star}] ||_{\infty} +u[\ell^{\star}]|_{\Xi}$, where $|| \cdot ||_{\infty}$ denotes in $L^{\infty}$ norm, $\delta=0.0001$ represents the noise level and $\xi_U$ is a vector of independent identically distributed uniform random variables on [-1,1]. The initial values are set as $\mu_1^0 = \cdots = \mu_{128}^0 = 0.35, \lambda_1^0 = \cdots = \lambda_{128}^0 = 1.5$. The iteration process continues for 482 steps until the termination criterion is met. Figure \ref{fig:J2} reflects the change in the value of $J$, which converges gradually with the iteration step number $k$. The numerically reconstructed Lam\'{e} coefficients $\mu_{h'}^{482}$ and $\lambda_{h'}^{482}$ along the $\phi$-direction is represented by the blue broken line in Figure \ref{one-lam2-mu2}, and the true values are indicated by the red dotted lines. At the 482nd iteration, the error results are: $\|\mu_{h'}^{482} - \mu^{\star}\|_{1}/\|\mu^{\star}\|_{1} = \text{8.223e-3}$ and $\|\lambda_{h'}^{482} - \lambda^{\star}\|_{1}/\|\lambda^{\star}\|_{1} = \text{4.132e-3}$.
\begin{figure}[H]
\centering
\includegraphics[scale=0.4]{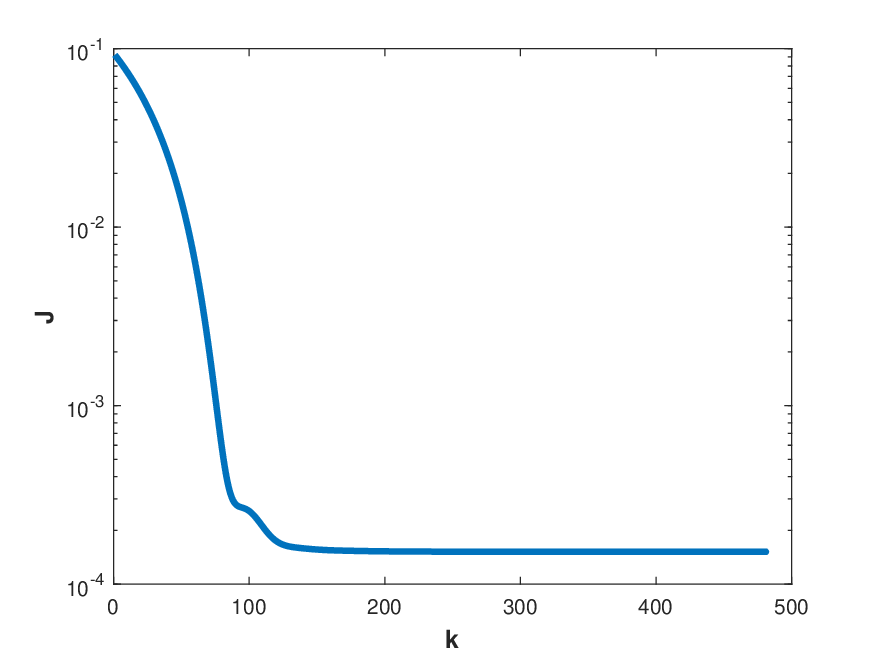}
\caption{Convergence of the objective function $J$ in Example \ref{example2}} 
\label{fig:J2}
\end{figure}
\begin{figure}[H]
\centering
\begin{subfigure}[b]{0.4\textwidth}
\includegraphics[width=\textwidth]{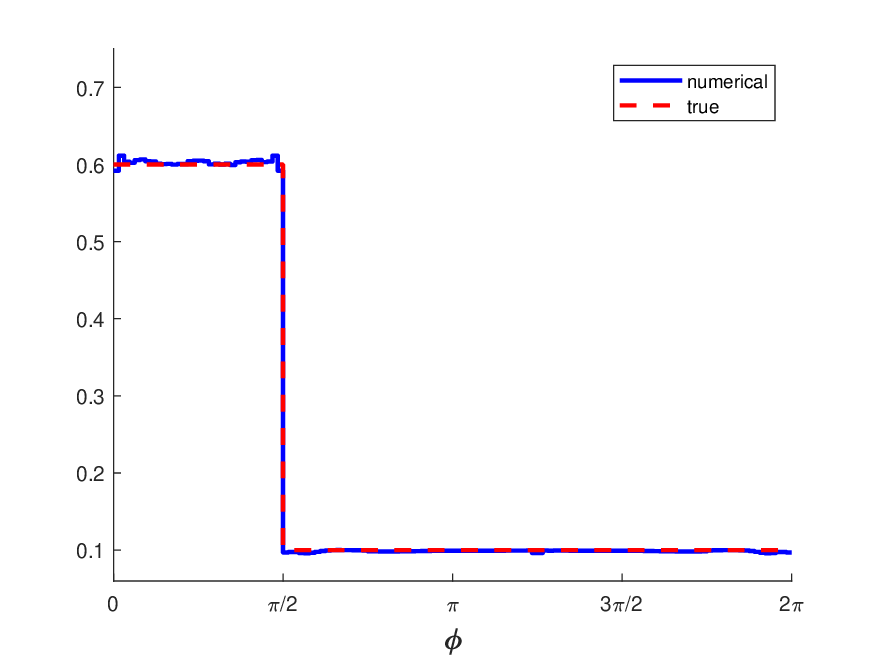}
\caption{$\mu$}
\end{subfigure}
\begin{subfigure}[b]{0.4\textwidth}
\includegraphics[width=\textwidth]{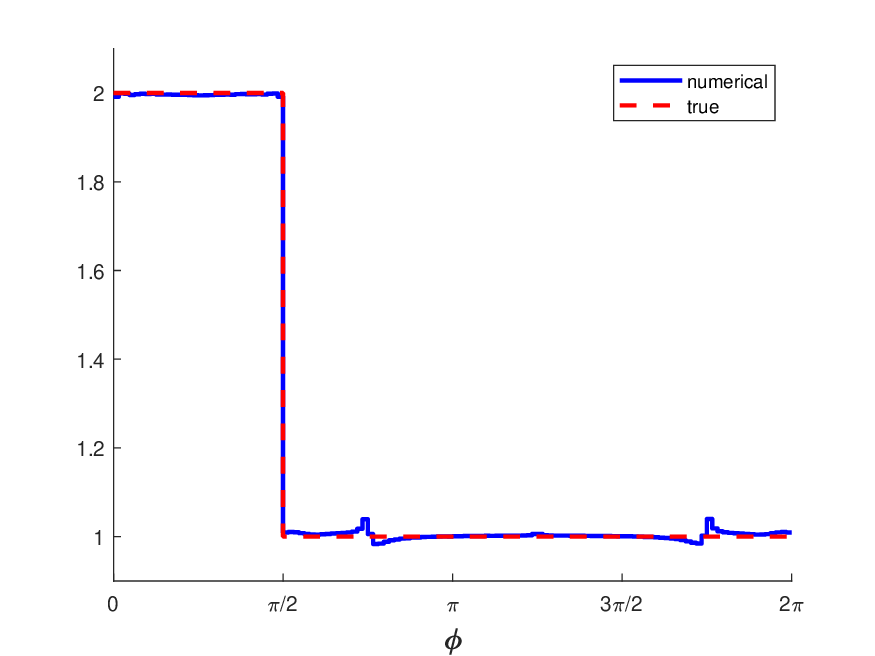}
\caption{$\lambda$}
\end{subfigure}
\caption{Reconstructed Lam\'e coefficients and true Lam\'e coefficients}
\label{one-lam2-mu2}
\end{figure}
\begin{example}\label{example3}
\end{example}
Building upon the forward problem analysis in Example \ref{ex3}, we now investigate the corresponding inverse elasticity problem to estimate the Lam\'{e} coefficients from displacement measurements.
The noisy measurement is defined as $z|_{\Xi} =\delta \xi_U || u[\ell^{\star}] ||_{\infty} +u[\ell^{\star}]|_{\Xi}$, where $|| \cdot ||_{\infty}$ denotes in $L^{\infty}$ norm, $\delta=0.0001$ represents the noise level and $\xi_U$ is a vector of independent identically distributed uniform random variables on [-1,1]. The initial values are set as $\mu_1^0 = \cdots = \mu_{128}^0 = 0.35, \lambda_1^0 = \cdots = \lambda_{128}^0 = 1.25$. The iteration process continues for 44 steps until the termination criterion is met. Figure \ref{fig:J3} reflects the change in the value of $J$, which converges gradually with the iteration step number $k$. The numerically reconstructed Lam\'{e} coefficients $\mu_{h'}^{44}$ and $\lambda_{h'}^{44}$ along the $\phi$-direction is represented by the blue broken line in Figure \ref{one-lam3-mu3}, and the true values are indicated by the red dotted lines. At the 44nd iteration, the error results are: $\|\mu_{h'}^{44} - \mu^{\star}\|_{1}/\|\mu^{\star}\|_{1} = \text{5.087e-3}$ and $\|\lambda_{h'}^{44} - \lambda^{\star}\|_{1}/\|\lambda^{\star}\|_{1} = \text{7.737e-3}$.

\begin{figure}[H]
\centering
\includegraphics[scale=0.4]{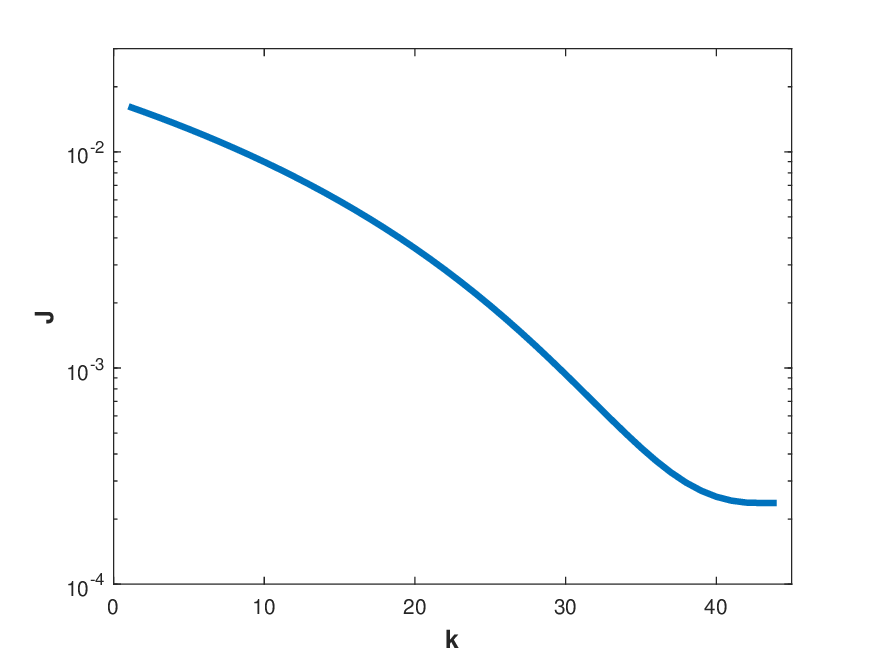}
\caption{Convergence of the objective function $J$ in Example \ref{example3}} %
\label{fig:J3}
\end{figure}
\begin{figure}[H]
\centering
\begin{subfigure}[b]{0.4\textwidth}
\includegraphics[width=\textwidth]{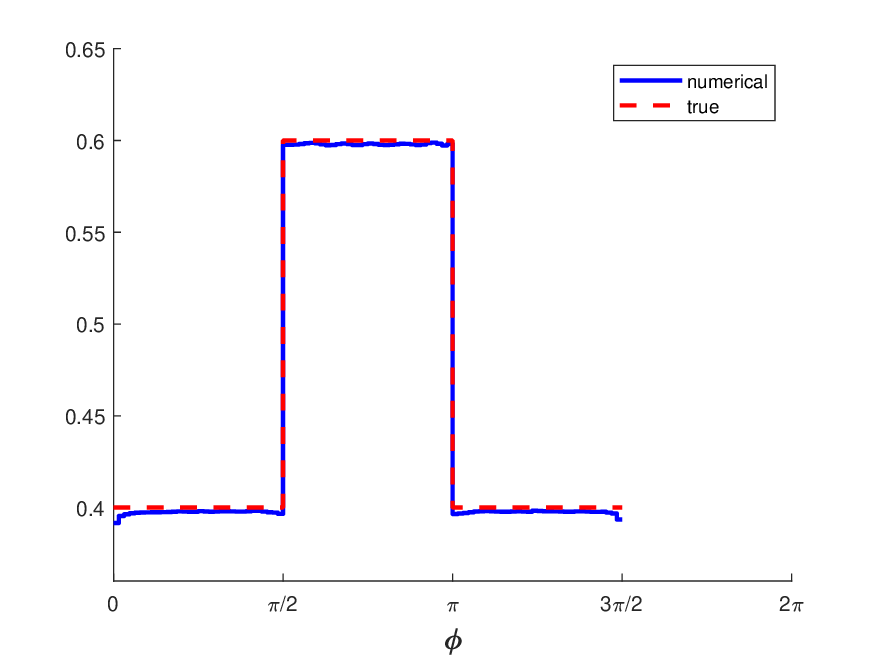}
\caption{$\mu$}
\end{subfigure}
\begin{subfigure}[b]{0.4\textwidth}
\includegraphics[width=\textwidth]{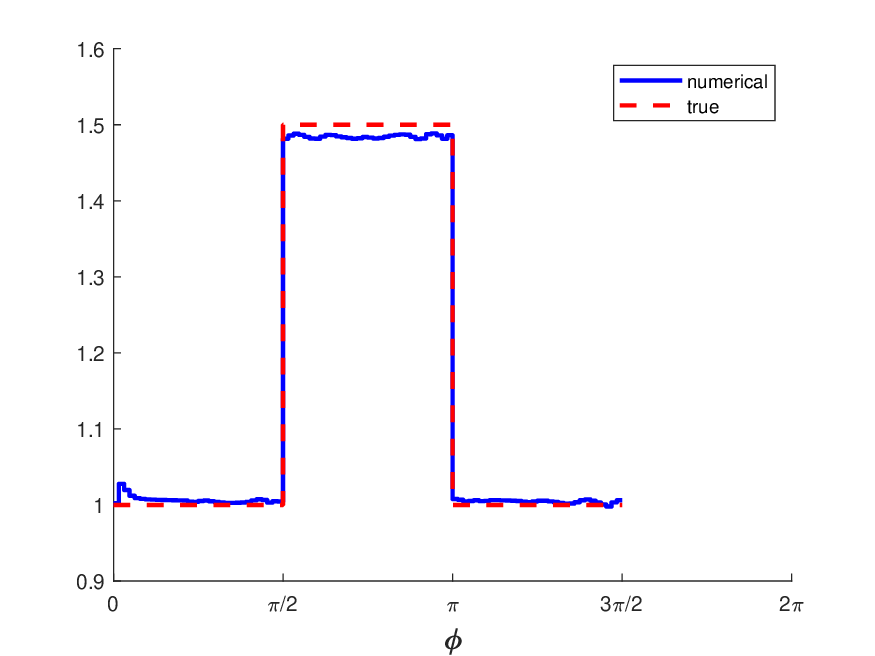}
\caption{$\lambda$}
\end{subfigure}
\caption{Reconstructed Lam\'e coefficients and true Lam\'e coefficients}
\label{one-lam3-mu3}
\end{figure}

\begin{example}\label{example4}
\end{example}
According to Example \ref{example3}, we will focus on the impact of measurement noise on the effectiveness and accuracy of the proposed Algorithm. We choose $\delta = 0.0003, 0.0007, 0.001, 0.01$. 
\begin{table}[H]
\centering
\caption{Error results under different noises}
\label{tab:case1}
\begin{tabular}{cccc}
\toprule
$\delta$ & $\tilde{k}$ & $\|\mu_{h'}^{\tilde{k}} - \mu^{\star}\|_{1} / \|\mu^{\star}\|_{1}$ & $\|\lambda_{h'}^{\tilde{k}} - \lambda^{\star}\|_1 / \|\lambda^{\star}\|_{1}$ \\
\midrule
0.0003 & 44 & 5.095e-3 & 7.310e-3 \\
0.0007 & 44 & 5.108e-3 & 7.355e-3 \\
0.001 & 44 & 5.509e-3 & 7.311e-3 \\
0.01 & 45 & 8.091e-3 & 1.702e-2 \\
\bottomrule
\end{tabular}
\end{table}
\begin{figure}[H]
\centering
\begin{subfigure}[b]{0.4\textwidth}
\includegraphics[width=\textwidth]{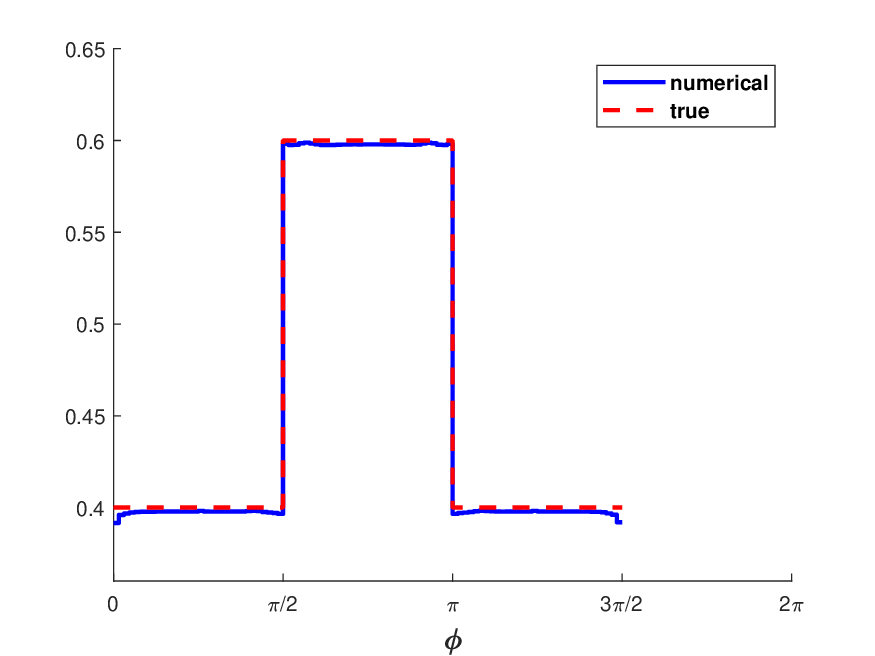}
\caption{$\mu$}
\end{subfigure}
\begin{subfigure}[b]{0.4\textwidth}
\includegraphics[width=\textwidth]{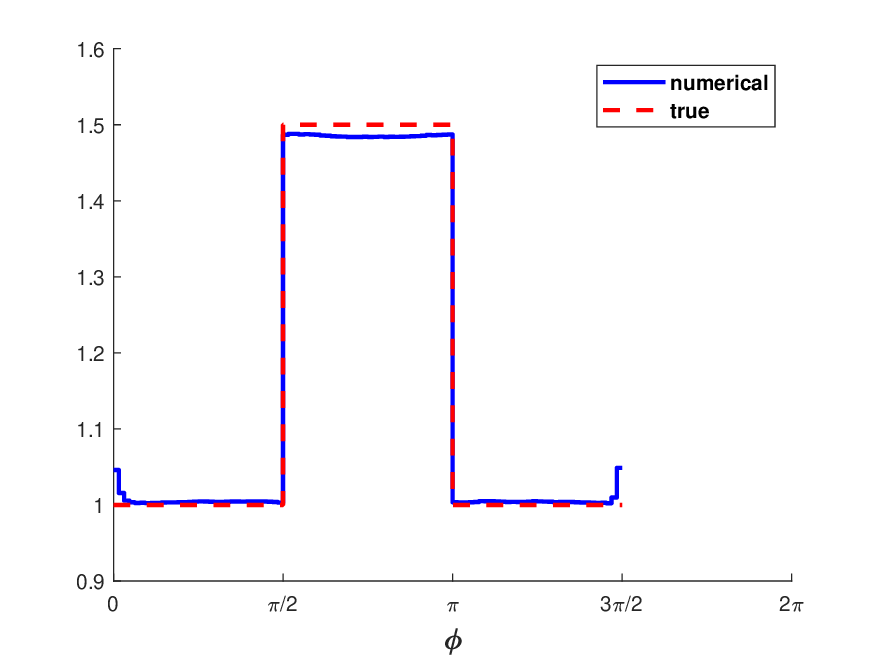}
\caption{$\lambda$}
\end{subfigure}
\caption{Reconstructed Lam\'e coefficients and true Lam\'e coefficients} \label{fig:four-lam-mu}
\label{case-one-lam-mu}
\end{figure}

For the case with noise level $\delta = 0.0001$, the reconstructed Lam\'e coefficients (see Figure \ref{one-lam3-mu3}) show excellent agreement with the true values. The corresponding results for higher noise levels $(\delta=0.0003,\delta=0.0007,\delta=0.001$ and $\delta=0.01)$ exhibited comparable accuracy in interface identification, though with slightly increased oscillations in the reconstructed parameter profiles near the material boundaries. (see Figure \ref{case-one-lam-mu}, here, only the image of $\delta=0.0003$ is shown). Table \ref{tab:case1} presents the convergence behavior under different noise levels $\delta$. In particular, our method is feasible and applicable to slight noise levels and different types of noise.
\section{Conclusion} \label{sec:conclusion}
This work proposes a novel computational framework that extends the direct-line method and domain decomposition techniques to linear elastic problems with multiple singular points in general domains and applies it to inverse elastic problems. Domain decomposition technology views the general domain as the union of two or more star-shaped subdomains. Assuming that the boundaries of each subdomain can be described by explicit $C^1$ parameter curves in polar coordinates, we transform the irregular star-shaped subdomain into a regular semi-infinite strip through coordinate transformation, thereby transforming our problem into a variational differential problem. And after discretizing the angular variable to obtain a semi-discrete approximation, we use the direct-line methods for analytical solution. Numerical results demonstrate that the eigenvalues converge rapidly to those of the exact elliptic operator, thus naturally capturing singularities, and this method possesses optimal error estimation capabilities. Numerical results confirm the accuracy and computational efficiency of the proposed method for general domains containing multiple singular points.

For inverse problems, this work addresses the inverse elasticity problem of composite materials by estimating piecewise constant Lam\'{e} coefficients through the minimization of an energy functional combined with total variation regularization. The optimization is implemented using the Adam algorithm. At every optimization step, the forward simulation is formulated as a linear elastic problem with multiple singular points on general domain. To handle geometric irregularities and potential multiple singularities, the direct-line method is employed as the forward solver. Numerical experiments demonstrate that the proposed approach accurately reconstructs both interface locations and Lam\'{e} coefficients distributions across all subdomains.


\section*{Acknowledgments}
This work was partially supported by the NSFC Project No. 12025104, Liaoning Provincial Natural Science Foundation Project No. 2025-BS-0410 and Liaoning Provincial Department of Education Project No. LJ212410147018.

\end{document}